\documentclass[12pt,reqno]{amsart}
\usepackage{amscd}
\usepackage{amssymb}
\usepackage{epsfig}
\usepackage{float}
\usepackage{url}
\usepackage[all]{xy}
\usepackage{graphicx}

\newcommand{\N}{{\mathbb N}}
\newcommand{\Z}{{\mathbb Z}}
\newcommand{\Q}{{\mathbb Q}}
\newcommand{\C}{{\mathbb C}}
\newcommand{\R}{{\mathbb R}}

\newcommand{\AAA}{{\mathcal A}}
\newcommand{\BB}{{\mathcal B}}
\newcommand{\CC}{{\mathcal C}}

\newcommand{\www}{\widetilde}

\newcommand{\oooo}{\overline}
\newcommand{\uuuu}{\underline}

\newcommand{\paa}{\partial}
\newcommand{\mmod}{{\rm mod\,}}
\newcommand{\Lo}{{\rm Lo}}
\newcommand{\Or}{{\rm Or}}
\newcommand{\sign}{{\rm sign}}

\DeclareMathOperator{\id}{id}
\DeclareMathOperator{\Imm}{Im}

\DeclareMathOperator{\pr}{pr}

\begin{document}

\theoremstyle{plain}
\newtheorem{lemma}{Lemma}[section]
\newtheorem{definition/lemma}[lemma]{Definition/Lemma}
\newtheorem{theorem}[lemma]{Theorem}
\newtheorem{proposition}[lemma]{Proposition}
\newtheorem{corollary}[lemma]{Corollary}
\newtheorem{conjecture}[lemma]{Conjecture}
\newtheorem{conjectures}[lemma]{Conjectures}

\theoremstyle{definition}
\newtheorem{definition}[lemma]{Definition}
\newtheorem{withouttitle}[lemma]{}
\newtheorem{remark}[lemma]{Remark}
\newtheorem{remarks}[lemma]{Remarks}
\newtheorem{example}[lemma]{Example}
\newtheorem{examples}[lemma]{Examples}
\newtheorem{notation}[lemma]{Notation}
\newtheorem{notations}[lemma]{Notations}

\title
[Integral monodromy of the cycle type singularities]
{The integral monodromy of the cycle type singularities}

\author{Claus Hertling \and Makiko Mase}

\address{Claus Hertling\\Universit\"at Mannheim\\ 
Lehrstuhl f\"ur algebraische Geometrie\\
B6, 26\\
68159 Mannheim, Germany}
\email{hertling@math.uni-mannheim.de}

\address{Makiko Mase\\Universit\"at Mannheim\\ 
Lehrstuhl f\"ur algbraische Geometrie\\
B6, 26\\
68159 Mannheim, Germany}
\email{mmase@mail.uni-mannheim.de}

\thanks{This work was funded by the Deutsche 
Forschungsgemeinschaft (DFG, German Research Foundation) 
-- 242588615}



\date{September 14, 2020}

\begin{abstract}
The middle homology of the Milnor fiber of a quasihomogeneous
polynomial with an isolated singularity is a $\Z$-lattice
and comes equipped with an automorphism of finite order,
the integral monodromy. 
Orlik (1972) made a precise conjecture, which would determine
this monodromy in terms of the weights of the polynomial.
Here we prove this conjecture for the cycle type singularities.
A paper of Cooper (1982) with the same aim contained 
two mistakes. Still it is very useful. 
We build on it and correct the mistakes. We give additional
algebraic and combinatorial results. 
\end{abstract}

\maketitle

\tableofcontents

\setcounter{section}{0}

\section{Introduction and main result}\label{c1}
\setcounter{equation}{0}

\noindent
The main objects of the paper are a quasihomogeneous 
singularity, the monodromy on its Milnor lattice, 
Orlik's conjecture for this monodromy, and our proof
in the case of a cycle type singularity. 
In order to make precise statements, we start with
some algebraic definitions.

\begin{definition}\label{t1.1}
%
(a) Start with a product $p\in\Z[t]$ of cyclotomic polynomials
which has only simple zeros. The {\it Orlik block} $\Or(p)$ 
is a pair $(H,h)$ where $H$ is a $\Z$-lattice of rank $\deg(p)$
and $h:H\to H$ is an automorphism of finite order 
with characteristic polynomial $p$ 
such that an element $a_0\in H$ with
\begin{eqnarray}\label{1.1}
H=\bigoplus_{j=0}^{\deg(p)-1}\Z\cdot h^j(a_0).
\end{eqnarray}
exists. Such an element is called a {\it generating element}.
The Orlik block $\Or(p)$ is up to isomorphism uniquely 
determined by $p$, which justifies the notion $\Or(p)$. 

\medskip 
(b) Consider a pair $(H,h)$ where $H$ is a $\Z$-lattice 
and $h:H\to H$ is an automorphism of finite order. 
It admits a {\it decomposition into Orlik blocks}
if it is isomorphic to a direct sum of Orlik blocks.

\medskip
(c) Consider a pair $(H,h)$ as in (b). 
Then the characteristic polynomial
$p_{H,h}$ of $h$ is a product of cyclotomic polynomials.
It has a unique decomposition $p_{H,h}=\prod_{i=1}^lp_i$
with $p_l|p_{l-1}|...|p_2|p_1$ and $p_l\neq 1$ and
all $p_i$ unitary and such that $p_1$ has only simple zeros. 
The pair $(H,h)$ admits a 
{\it standard decomposition into Orlik blocks},
if an isomorphism $(H,h)\cong \bigoplus_{i=1}^l\Or(p_i)$ exists.
\end{definition}

A polynomial $f\in\C[x_1,...,x_n]$ is called quasihomogeneous if
for some weight system $(w_1,...,w_n)$ with 
$w_i\in(0,1)\cap\Q$ each monomial in $f$ has weighted degree 1. 
It is called an isolated
quasihomogeneous singularity if it is quasihomogeneous and
the functions 
$\frac{\paa f}{\paa x_1},...,\frac{\paa f}{\paa x_n}$
vanish simultaneously only at $0\in\C^n$. 
Then the Milnor lattice $H_{Mil}:=H_{n-1}^{(red)}(f^{-1}(1),\Z)$
(here $H_{n-1}^{(red)}$ means the reduced homology in the 
case $n=1$ and the usual homology in the cases $n\geq 2$)
is a $\Z$-lattice of some finite rank $\mu\in\N$, which is called
the {\it Milnor number} \cite{Mi68}. It comes equipped with
an automorphism $h_{Mil}$ of finite order, the {\it monodromy}.

Orlik conjectured the following.

\begin{conjecture}\label{t1.2}
(Orlik's conjecture \cite[Conjecture 3.1]{Or72})
For any isolated quasihomogeneous singularity, the pair
$(H_{Mil},h_{Mil})$ admits a standard decomposition into Orlik blocks.
\end{conjecture}

Here we will prove this conjecture for the cycle type 
singularities. A cycle type singularity is a polynomial
$f\in\C[x_1,...,x_n]$ in $n\geq 2$ variables of the following
shape,
\begin{eqnarray}\label{1.2} 
f=x_1^{a_1}x_2+x_2^{a_2}x_3+...+x_{n-1}^{a_{n-1}}x_n
+x_n^{a_n}x_1 \quad
\textup{where }a_1,...,a_n\in\N,
\end{eqnarray}
($\N=\{1,2,3,...\}$) and where 
\begin{eqnarray}
\textup{for even }n
\textup{ neither }a_j=1
\textup{ for all even }j\textup{ nor }a_j=1
\textup{ for all odd }j.\label{1.3}
\end{eqnarray}
It is an isolated quasihomogeneous singularity. 
Define $d:=\prod_{j=1}^na_j-(-1)^n$. Then 
(see e.g. Lemma 4.1 in \cite{HZ19}) 
\begin{eqnarray*}
\mu=\prod_{j=1}^n a_j =d+(-1)^n,
\end{eqnarray*}
and there are natural numbers $v_1,...,v_n$ (which are given
in \eqref{t5.1}) such that 
$(w_1,...,w_n)=(\frac{v_1}{d},...,\frac{v_n}{d})$ is the
unique weight system for which $f$ is quasihomogeneous 
of weighted degree 1. They satisfy also 
$\gcd(v_1,d)=...=\gcd(v_n,d)$. Define 
$b:=d/\gcd(v_1,d)\in\N$. An easy calculation which builds on
the formula in \cite{MO70} for the characteristic polynomial
in terms of $(w_1,...,w_n)$ (see e.g. Lemma 4.1 in \cite{HZ19})
shows that the characteristic polynomial of $h_{Mil}$
on $H_{Mil}$ is
\begin{eqnarray}\label{1.4}
p_{H_{Mil},h_{Mil}}=(t^b-1)^{\gcd(d,v_1)}\cdot (t-1)^{(-1)^n}.
\end{eqnarray}
Therefore Orlik's conjecture says here the following.

\begin{theorem}\label{t1.3}
For a cycle type singularity as above,
\begin{eqnarray*}\label{1.5}
(H_{Mil},h_{Mil})\cong \left\{\begin{array}{ll} 
(\gcd(d,v_1)-1)\Or(t^b-1)\oplus\Or(\frac{t^b-1}{t-1})
&\textup{if }n\textup{ is odd},\\
\gcd(d,v_1)\Or(t^b-1)\oplus\Or(t-1)
&\textup{if }n\textup{ is even}.\end{array}\right.
\end{eqnarray*}
\end{theorem}

Our proof in section \ref{c5} 
builds on Cooper's work \cite{Co82}. 
By \cite{Mi68}, the set
$\oooo{F_0}:=f^{-1}(\R_{\geq 0})\cap S^{2n-1}\subset S^{2n-1}
\subset\C^n$
is diffeomorphic to the Milnor fiber $f^{-1}(1)$.
Cooper studied for the cycle type singularities a beautiful subset
$G\subset\oooo{F_0}$, which is a probably
a deformation retract of the Milnor fiber (Cooper's Lemma 3
is slightly weaker). He considered certain {\it cells}
from which this set $G$ is built up and which allow to filter
$G$ by a sequence of subsets $G=G_n\supset G_{n-1}\supset ...
\supset G_1\supset G_0=\emptyset$. 
Finally, he studied a spectral sequence which comes from
this filtration.

Cooper claimed to have proved Orlik's conjecture for the
cycle type singularities. But his paper contains two serious
mistakes. The second one leads in the case of even $n$ to
the wrong claim
$$(H_{Mil},h_{Mil})\cong(\gcd(d,v_1)-1)\Or(t^b-1)
\oplus\Or(\frac{t^b-1}{t-1})\oplus 2\Or(t-1).$$
The right hand side is a decomposition into Orlik blocks, 
but not a standard decomposition.

We will discuss the relation of this paper to Cooper's work
and the two mistakes in the Remarks \ref{t5.1}.
His work is the basis for the sections \ref{c3} and \ref{c4}
below. Section \ref{c2} gives new algebraic results.
Section \ref{c3} gives the set $G$, its cells 
and an inductive construction of cycles 
of which only the beginning is in \cite{Co82}. 
Section \ref{c4} makes good use of the spectral 
sequence which Cooper considered and determines $H_{n-1}(G,\Z)$. 
A combination of the results of the sections \ref{c2},
\ref{c3} and \ref{c4} and a discussion of the monodromy
proves Theorem \ref{t1.3} in section \ref{c5}.

\section{Algebraic results}\label{c2}
\setcounter{equation}{0}

\noindent
The following two lemmata \ref{t2.3} and \ref{t2.4} 
are elementary. They will be used in the proof
of Orlik's conjecture for cycle type singularities with an
even number of variables. 
They say something about the $\Z$-lattices
$H^{(d,c)}$ with automorphisms $h^{(d,c)}$ of finite order,
which are defined in Definition \ref{t2.1}.

\begin{definition}\label{t2.1}
Let $d\in\N$ and $c\in\Z$. Define the pair
$\Lo^{(d,c)}=(H^{(d,c)},h^{(d,c)})$ as follows.
$H^{(d,c)}=\Z\cdot\gamma\oplus
\bigoplus_{j=1}^{d-1}\Z\cdot \delta_j$ is a $\Z$-lattice of
rank $d$. Define additionally $\delta_d:=c\cdot \gamma
-\sum_{j=1}^{d-1}\delta_j\in H^{(d,c)}$. 
Then $h^{(d,c)}:H\to H$ is the automorphism
of finite order $d$ which is defined by
\begin{eqnarray}\label{2.1}
h^{(d,c)}:\ \gamma\mapsto\gamma,\quad \delta_d
\mapsto \delta_1,\quad 
\delta_j\mapsto \delta_{j+1}\textup{ for }j\in\{1,...,d-1\}.
\end{eqnarray}
\end{definition}

\begin{remark}\label{t2.2}
The characteristic polynomial of $h^{(d,c)}$ is $t^d-1$.
If $c\neq 0$, then 
$\sum_{j=1}^d\Z\cdot\delta_j=\bigoplus_{j=1}^d\Z\cdot\delta_j$
is an $h^{(d,c)}$-invariant sublattice of index $|c|$ in 
$H^{(d,c)}$, and 
$(\bigoplus_{j=1}^d\Z\cdot\delta_j,h^{(d,c)})
\cong\Or(t^d-1)$.
Therefore $\Lo(d,1)\cong\Or(t^d-1)$. 
If $c=0$ then the summands $\Z\cdot\gamma$ and
$\bigoplus_{j=1}^{d-1}\Z\cdot\delta_j$ of $H^{(d,c)}$ 
are $h^{(d,c)}$-invariant
with $(\Z\cdot\gamma,h^{(d,c)})\cong \Or(t-1)$ and
$(\bigoplus_{j=1}^{d-1}\Z\cdot\delta_j,h^{(d,c)})
\cong \Or(\frac{t^d-1}{t-1})$. 
\end{remark}

\begin{lemma}\label{t2.3}
Let $d,v\in\N$, $c\in\Z$ and $b:=d/\gcd(d,v)\in\N$. Then
\begin{eqnarray}\label{2.2}
(H^{(d,c)},(h^{(d,c)})^v)\cong
(\gcd(d,v)-1)\cdot\Or(t^b-1)\oplus \Lo^{(b,c)}.
\end{eqnarray}
\end{lemma}

{\bf Proof:} Write $h:=(h^{(d,c)})^v$. 
The elements $\delta_1,...,\delta_d$ can be renumbered
to elements $\www\delta_1,...,\www\delta_d$
(i.e. $\{\delta_1,...,\delta_d\}=
\{\www\delta_1,...,\www\delta_d\}$) such that these form
$\gcd(d,v)$ many {\it cycles of length $b$} with respect to $h$:
\begin{eqnarray*}
h: && \www\delta_{ab+1}\mapsto \www\delta_{ab+2}\mapsto 
... \mapsto \www\delta_{ab+b}\mapsto \www\delta_{ab+1}\\
&& \textup{ for }a\in\{0,1,...,\gcd(d,v)-1\}.
\end{eqnarray*}
Define
\begin{eqnarray*}
\beta_j:=\sum_{a=0}^{\gcd(d,v)-1}\www\delta_{ab+j}
\quad\textup{ for }j\in\{1,...,b\}.
\end{eqnarray*}
Then these also form a cycle of length $b$ with respect to $h$,
and their sum is $c\gamma$, 
\begin{eqnarray*}
h:\ \beta_1\mapsto \beta_2\mapsto...\mapsto \beta_b\mapsto
\beta_1,\\
\sum_{j=1}^b\beta_j=\sum_{j=1}^d\www\delta_j
=\sum_{j=1}^d\delta_j=c\gamma.
\end{eqnarray*}
We obtain
\begin{eqnarray*}
(H^{(d,c)},h)&\cong& 
\bigoplus_{a=0}^{\gcd(d,v)-2}\left(\bigoplus_{j=1}^b\Z\cdot
\www\delta_{ab+j},h\right)\oplus
\left(\Z\cdot\gamma+\sum_{j=1}^b\Z\cdot\beta_j,h\right)\\
&\cong& 
(\gcd(d,v)-1)\Or(t^b-1)\oplus \Lo^{(b,c)}. \hspace*{3cm}\Box
\end{eqnarray*}

\begin{lemma}\label{t2.4}
Let $d\in\N$ and $c,\www c\in\Z$. The following three conditions
are equivalent.
\begin{list}{}{}
\item[(i)]
$\Lo^{(d,c)}\cong \Lo^{(d,\www c)}.$
\item[(ii)]
$\Lo^{(d,c)}\oplus \Or(t-1)\cong \Lo^{(d,\www c)}\oplus \Or(t-1).$
\item[(iii)]
$\gcd(d,c)=\gcd(d,\www c).$
\end{list}
And
\begin{eqnarray}\label{2.3}
&&\Lo^{(d,c)}\oplus\Or(t-1)\cong \Or(t^d-1)\oplus\Or(t-1)\\
&\iff& \gcd(d,c)=1.\nonumber
\end{eqnarray}
\end{lemma}

{\bf Proof:} 
Keep the notations $\delta_1,...,\delta_d,\gamma$ 
of Definition \ref{t2.1} for the elements of $H^{(d,c)}$.
And extend them by $\delta_j:=\delta_{j_0}$ if
$j\in\Z-\{1,...,d\}$, $j_0\in\{1,...,d\}$ and $d|(j-j_0)$. 

\medskip
(iii)$\Rightarrow$(i): 
We start with $\Lo^{(d,c)}$ for some $c\in\N$ with $c|d$.
We will present a construction which leads to certain 
$\www c\in\Z$ with $\Lo^{(d,c)}\cong\Lo^{(d,\www c)}$. 
Then we will show that these are all $\www c\in\Z$ with
$\gcd(d,\www c)=c$. 

Choose $a\in\N$ with $\gcd(a,d)=1$, and choose $b\in\Z$. 
Define
\begin{eqnarray*}
\www\gamma:=\gamma,\quad 
\www\delta_j:=b\gamma+\sum_{i=0}^{a-1}\delta_{j+i}\quad
\textup{for }j\in\Z,
\end{eqnarray*}
so that $\www\delta_{j_1}=\www\delta_{j_2}$ if $d|(j_1-j_2)$. 
Of course, $h^{(d,c)}$ acts by 
\begin{eqnarray}\label{2.4}
h^{(d,c)}:\ \www\gamma\mapsto\www\gamma,\quad\www\delta_j\mapsto
\www\delta_{j+1},
\end{eqnarray}
and we have 
\begin{eqnarray}\label{2.5}
\sum_{j=1}^d\www\delta_j=(bd+ac)\cdot\www\gamma.
\end{eqnarray}
Furthermore, the condition $\gcd(a,d)=1$ implies
\begin{eqnarray}\label{2.6}
\Z\cdot\www\gamma+\sum_{j=1}^d\Z\cdot\www\delta_j=
H^{(d,c)} 
\left(=\Z\cdot\gamma+\sum_{j=1}^d\Z\cdot\delta_j\right).
\end{eqnarray}
To see this, choose $b_1,b_2\in\N$ with $ab_1-db_2=1$.
Then for $j\in\Z$
\begin{eqnarray*}
-(b_1b+b_2c)\www\gamma+\sum_{k=0}^{b_1-1}\www\delta_{j+ak}
=-b_1b\www\gamma-\sum_{k=1}^{db_2}\delta_{j+k}
+b_1b\www\gamma+\sum_{k=0}^{db_2}\delta_{j+k}
=\delta_j,
\end{eqnarray*}
which shows \eqref{2.6}. Together, \eqref{2.4}, \eqref{2.5}
and \eqref{2.6} give
\begin{eqnarray}\label{2.7}
\Lo^{(d,c)}&\cong&\Lo^{(d,bd+ac)}.
\end{eqnarray}

Now choose any $\www c\in\Z$ with $\gcd(d,\www c)=c$.
It remains to see that there exist $a\in\N$ and $b\in\Z$ with
$\gcd(a,d)=1$ and $\www c=bd+ac$. 

For any integer $m\in\Z-\{0\}$, write 
$m=\frac{m}{|m|}\cdot
\prod_{p\textup{ prime number}}p^{v_p(m)}$, where
$v_p(m)\in\Z_{\geq 0}$. We choose 
\begin{eqnarray*}
\www{b}&:=& \prod_{p\textup{ prime number with }
v_p(d)=v_p(\www c)>0}p,\\
b&:=& -\www{b}-|\www{c}|\in\Z_{<0},\\
a&:=&\frac{\www c}{c}-b\cdot\frac{d}{c}
=\frac{\www c}{c}+\www{b}\frac{d}{c}+
\frac{|\www{c}|}{c}d\in\N.
\end{eqnarray*}
Then $\www c=bd+ac$. For a prime number $p$ with
$v_p(d)>v_p(\www c)$, $v_p(\frac{\www c}{c})=0$ and
$v_p(\frac{d}{c})>0$ and $v_p(a)=0$. 
For a prime number $p$ with $v_p(d)=v_p(\www c)>0$
$v_p(\frac{\www c}{c})=0$ and 
$v_p(b)>0$ and $v_p(a)=0$. 
For a prime number $p$ with $0<v_p(d)<v_p(\www c)$,
$v_p(\frac{\www{c}}{c})>0$ and $v_p(b\frac{d}{c})=0$ and
$v_p(a)=0$. Therefore $\gcd(a,d)=1$.

\medskip
(i)$\Rightarrow$(ii): This is trivial.

\medskip
(ii)$\Rightarrow$(iii): We proved already (iii)$\Rightarrow$(i)
and (i)$\Rightarrow$(ii). Therefore, in (ii) we can suppose
$c=\gcd(d,c)$ and $\www c=\gcd(d,\www c)$. 
Then we have to show $c=\www c$.

Write the elements in Definition \ref{t2.1} 
for $H^{(d,\www c)}$ with a tilde, so as $\www\delta_1,...,
\www\delta_d,\www\gamma$. 
Write generators of $\Or(t-1)$ on the left hand side
respectively right hand side of (ii) as 
$\beta$ respectively $\www\beta$. The automorphisms
of the left hand side respectively right hand side of
(ii) which extend $h^{(d,c)}$ respectively
$h^{(d,\www c)}$ by $\id$ on $\Or(t-1)$, are called
$h$ respectively $\www h$. 

Let 
$$g:\Lo^{(d,c)}\oplus\Or(t-1)\to
\Lo^{(d,\www c)}\oplus\Or(t-1)$$ 
be an isomorphism. Then 
\begin{eqnarray}\label{2.8}
\www h\circ g=g\circ h.
\end{eqnarray}
This and $\ker(h-\id)=\Z\cdot\gamma\oplus\Z\cdot\beta$
and 
$\ker(\www h-\id)=\Z\cdot\www\gamma\oplus\Z\cdot\www\beta$ imply
\begin{eqnarray*}
g(\gamma)=b_1\www\gamma+b_2\www\beta,\quad
g(\beta)=b_3\www\gamma+b_4\www\beta\quad
\textup{with }\begin{pmatrix}b_1&b_2\\b_3&b_4\end{pmatrix}
\in GL(2,\Z).
\end{eqnarray*}
Because of \eqref{2.8} and the action of $h$ on
$\delta_1,...,\delta_d$ and $\www{h}$ on $\www{\gamma}$,
$\www{\beta}$, $\www{\delta}_1,...,\www{\delta}_d$, 
numbers $a_1,...,a_{d+2}\in\Z$ with
\begin{eqnarray*}
g(\delta_j)=a_{d+1}\www\gamma+a_{d+2}\www\beta
+\sum_{i=1}^d a_i\www\delta_{i+j-1}
\quad\textup{for }j\in\{1,...,d\}
\end{eqnarray*}
exist. Write $a:=\sum_{i=1}^da_i$. Then
\begin{eqnarray*}
c(b_1\www\gamma+b_2\www\beta) &=& g(c\gamma)
=g(\sum_{j=1}^d\delta_j)\\
&=&da_{d+1}\www\gamma+da_{d+2}\www\beta+
\left(\sum_{i=1}^da_i\right)
\cdot\sum_{j=1}^d\www\delta_j \\
&=&(da_{d+1}+a\www c)\www\gamma + da_{d+2}\www\beta,\\
\textup{so }\quad 
cb_1&=& da_{d+1}+a\www c,\quad cb_2=da_{d+2}.
\end{eqnarray*}
$\pm 1= b_1b_4-b_2b_3$ implies $\gcd(b_1,b_2)=1$.
Therefore $c$ is in the ideal generated by $d$ and $\www c$,
which is the ideal generated by $\www c$, because $\www c|d$. 
Thus $\www c|c$.
Because of the symmetry of the situation also $c|\www c$. 
We obtain $c=\www c$. This shows (ii)$\Rightarrow$(iii).

\medskip
Now the equivalence of (i), (ii) and (iii) is proved.
It rests to see \eqref{2.3}. It follows with the isomorphy 
$\Lo^{(d,1)}\cong\Or(t^d-1)$ in Remark \ref{t2.2}
and the equivalence of (i), (ii) and (iii).
\hfill$\Box$

\section{Cells, chains and cycles}\label{c3}
\setcounter{equation}{0}

\noindent
Throughout this section we fix 
a cycle type singularity $f(x_1,...,x_n)$ as in 
\eqref{1.2} with $n\geq 2$ and $a_1,...,a_n\in\N$ 
with \eqref{1.3}. By \cite{Mi68}, the Milnor fiber 
$f^{-1}(1)\subset\C^n$
and the set $\oooo{F_0}:= f^{-1}(\R_{\geq 0})\cap S^{2n-1}
\subset S^{2n-1} \subset\C^n$ are diffeomorphic. 
Cooper \cite{Co82} considers the subset 
$G\subset \oooo{F_0}$
which is defined as follows,
\begin{eqnarray}\label{3.1}
G=\{z\in S^{2n-1}\,|\, 
z_j^{a_j}z_{j+1}\geq 0\ \forall\ j\in\{1,...,n-1\},\ 
z_n^{a_n}z_1\geq 0\}.
\end{eqnarray}
He conjectures that $G$ is a deformation retract of 
$\oooo{F_0}$. He proves a slightly weaker deformation lemma
(stated at the end of 3. in \cite{Co82}) 
which implies especially that the inclusion map 
$i_g:G\hookrightarrow \oooo{F_0}$ induces epimorphisms 
in homology. For him and for us, this property suffices. 

Cooper builds $G$ up from certain {\it cells}.
We will need these cells, and also refinements of them.
For this, quite some notations are needed. They are given now.

\begin{notations}\label{t3.1}
$N:=\{1,...,n\}$, so that $\R^N=\R^n$. The map
\begin{eqnarray*}
\uuuu{e}:\R^n\to T^n:=(S^1)^n\subset\C^n,\quad
(r_1,...,r_n)\mapsto  (e^{2\pi i r_1},...,e^{2\pi ir_n}),
\end{eqnarray*}
induces an isomorphism $\uuuu{e}_T:\R^n/\Z^n\to T^n$. 
Also the projection $\pr_T:\R^n\to\R^n/\Z^n$ will be useful.
Then $\uuuu{e}=\uuuu{e}_T\circ \pr_T$. The following 
binary operation $\odot$ is not standard, but it will also 
be useful,
\begin{eqnarray*}
\odot: \R_{\geq 0}^n\times(S^1)^n\to\C^n, 
(a_1,...,a_n)\odot (b_1,...,b_n):=(a_1b_1,...,a_nb_n).
\end{eqnarray*}

Given a finite tuple $(v_1,...,v_k)\in(\R^n)^k$ of vectors 
in $\R^n$, we consider the subset of $\R^n$
\begin{eqnarray*}
C(v_1,...,v_k):=\{\sum_{i=1}^kt_iv_i\,|\, t_1,...,t_k\in[0,1]\}.
\end{eqnarray*}
If $(v_1,...,v_k)$ is a tuple of linearly independent vectors
(which will be the case almost always), 
$C(v_1,...,v_k)$ is a hypercube of dimension $k$. 
Then it inherits an orientation from the ordered tuple
$(v_1,...,v_k)$. And then its boundary is
\begin{eqnarray}\label{3.2}
\partial C(v_1,...,v_k)&=& \sum_{i=1}^k
(-1)^{i-1}C(v_1,...,\widehat{v_i},...,v_k)\\
&&-\sum_{i=1}^k (-1)^{i-1}(v_i+C(v_1,...,\widehat{v_i},...,v_k)),
\nonumber
\end{eqnarray}
where $\widehat{v_i}$ means that $v_i$ is erased in the 
tuple.  
The following observation will be useful, because we will
consider the images of hypercubes under $\pr_T$ in $\R^n/\Z^n$. 
If $v_1\in\Z^n$, then
\begin{eqnarray}\label{3.3}
\pr_T(C(v_1,...,v_k))&=&\pr_T(C(v_1,v_2+\lambda_2v_1,...,
v_k+\lambda_kv_1))\\
&&\textup{for any }\lambda_2,...,\lambda_k\in\R,\nonumber
\end{eqnarray}
and similarly for $v_j\in\Z^n$ instead of  $v_1\in\Z^n$.

The exponents $a_1,...,a_n\in\N$ in the monomials of the 
cycle type singularity are the source of useful integers 
and vectors of integers:
\begin{eqnarray*}
\textup{Recall }\quad \mu&=&a_1\cdot...\cdot a_n,\quad  
d:=\mu-(-1)^n.\\
\textup{For }k\in N:\quad \uuuu{a}_k&:=&(-1)^{k-1}a_1\cdot ...
\cdot a_{k-1}\in\Z,\quad\textup{ so }\uuuu{a}_1=1.\\
\uuuu{a}_{n+1}&:=& (-1)^na_1\cdot...\cdot a_n=(-1)^n\mu
=(-1)^nd+1.\\
\textup{For }k\in N:\quad \uuuu{b}_k&:=&
(\uuuu{a}_1,...,\uuuu{a}_{k-1},0,...,0)\in\Z^n,
\quad\textup{so }\uuuu{b}_1=(0,...,0).\\
\uuuu{b}_{n+1}&:=& (\uuuu{a}_1,\uuuu{a}_2,...,\uuuu{a}_n)\in\Z^n.
\\
\textup{For }k\in N:\quad \uuuu{c}_k&:=&
(\uuuu{a}_{n+1}\uuuu{a}_{1},...,\uuuu{a}_{n+1}\uuuu{a}_{k-1},
\uuuu{a}_k,\uuuu{a}_{k+1},...,\uuuu{a}_n)\\
&=&\uuuu{a}_{n+1}\cdot\uuuu{b}_k
+(0,...,0,\uuuu{a}_k,...,\uuuu{a}_n),\\
&=&(-1)^nd\cdot\uuuu{b}_k+\uuuu{c}_1\in\Z^n,
\quad\textup{especially }\uuuu{c}_1=\uuuu{b}_{n+1}.
\end{eqnarray*}
Also the following vectors in $\Q^n$ will be used,
\begin{eqnarray*}
\textup{for }j\in\{1,...,d\}:\quad
\uuuu{p}_j&:=& \frac{j}{d}\cdot \uuuu{c}_1\in\Q^n.
\end{eqnarray*}
In fact, their images $\pr_T(\uuuu{p}_j)$ in $\R^n/\Z^n$
will be the only 0-chains in $\R^n/\Z^n$ which we will need.  
The reason is the equality 
\begin{eqnarray}\label{3.4}
\pr_T(\uuuu{p}_{j+1})=\pr_T(\uuuu{p}_j+\frac{1}{d}\uuuu{c}_k)
\quad \textup{for any }k\in N.
\end{eqnarray}

For any subset $A\subset N$ with $A\neq\emptyset$, define 
\begin{eqnarray*}
\R^A&:=&\{r=(r_1,...,r_n)\in\R^n\,|\, r_j=0\textup{ for }j
\notin A\}\subset \R^n,\\
\Z^A&\subset& \Z^n\textup{ and }\R^A_{\geq 0}\subset\R^n_{\geq 0}
\textup{ and }\C^A\subset\C^n \quad\textup{analogously},\\
\Delta_A&:=& \R^A_{\geq 0}\cap S^{2n-1}\subset S^{2n-1},\\
\pr_A&:&\R^n\to\R^A\quad\textup{the projection}.
\end{eqnarray*}
$\Delta_A$ is a deformation of a simplex of dimension $|A|-1$.
We call it {\it a deformed simplex}. Its boundary consists of 
the deformed simplices $\Delta_B$ with $B\subsetneqq A$.

We want to consider  $N$ and its subsets $A$ as cyclic: 
$1$ follows $n$. In order to write this down, we denote
by $(k)_{\mmod n}\in N$ for $k\in\Z$ the number with
$n|(k-(k)_{\mmod n})$. 

For $A\subsetneqq N$ with $A\neq\emptyset$, 
the {\it blocks} are the maximal sequences 
$k,(k+1)_{\mmod n},...,(k+l)_{\mmod n}$ within $A$. 
Their number is called $b(A)$. And the {\it block beginnings} are
the first numbers in the blocks. Explicitly, they are the
numbers in the set
\begin{eqnarray*}
\{k_1^A,...,k_{b(A)}^A\}
=\{k\in A\,|\,(k-1)_{\mmod n}\notin A\},\\
\textup{with }1\leq k_1^A<...<k_{b(A)}^A\leq n.
\end{eqnarray*}
The {\it gaps} of $A$ are the blocks
of $N-A$. Also their number is $b(A)$. 
A set $A\subsetneqq N$ is {\it thick} if $b(A)=n-|A|$.
Equivalent is that each gap of $A$ consists of a single
number. 
For thick $A\subset N$  define the sign
\begin{eqnarray*}
\sign(A):= \sign\begin{pmatrix}1&...&|A|&|A|+1&...&n\\
\alpha_1&...&\alpha_{|A|}&k_1^A&...&k_{b(A)}^A\end{pmatrix}
\in\{\pm 1\}\quad\textup{where}\\
N-\{k_1^A,...,k_{b(A)}^A\}=\{\alpha_1,...,\alpha_{|A|}\}
\textup{ with } \alpha_1<...<\alpha_{|A|}.
\end{eqnarray*}

A set $A\subsetneqq N$ with $A\neq\emptyset$ is {\it almost thick}
if $b(A)=n-|A|-1$. Equivalent is that one gap consists of two
numbers and each other gap consists of a single number.
For an almost thick set $B$ with gap $\{k_0,(k_0+1)_{\mmod n}\}$
of two elements, denote
\begin{eqnarray*}
B^{(1)}:= B\cup\{k_0\},\quad B^{(2)}:=B\cup\{(k_0+1)_{\mmod n}\}.
\end{eqnarray*}
$B^{(1)}$ and $B^{(2)}$ are the unique thick sets with
$B\subset B^{(i)}$ and $|B^{(i)}|=|B|+1$. They satisfy
$b(B)=b(B^{(i)})$. 

For $A\subsetneqq N$ with $A\neq\emptyset$, we will define
a subtorus $T_A\subset T^n$ below. For this we define for each
block beginning $k_j^A$ in $A$ a vector of integers,
\begin{eqnarray*}
\uuuu{d}_j^A&:=& \uuuu{a}_{k_j^A}^{-1}(\uuuu{b}_{k_{j+1}^A}
-\uuuu{b}_{k_j^A})\quad\textup{for }j\in\{1,...,b(B)-1\}\\
&=&(0,...,0,1,(-a_{k_j^A}),...,
(-a_{k_j^A})...(-a_{k_{j+1}^A-2}),0,...,0)\in\Z^n,\\
\uuuu{d}_{b(A)}^A&:=& \uuuu{a}_{k_{b(A)}^A}^{-1}
(\uuuu{a}_{n+1}\uuuu{b}_{k_1^A}+\uuuu{c}_1
-\uuuu{b}_{k_{b(A)}^A})\\
&=&(...,\uuuu{a}_{n+1}\uuuu{a}_{k_{b(A)}^A}^{-1}
\uuuu{a}_{k_1^A-1},0,...,0,1,(-a_{k_{b(A)}^A}),...,
\uuuu{a}_n\uuuu{a}_{k_{b(A)}^A}^{-1})\in\Z^n.
\end{eqnarray*}
Then the subtorus $T_A\subset T^n$ is the set
\begin{eqnarray*}
T_A&:=& \uuuu{e}(C(\uuuu{d}_1^A,...,\uuuu{d}_{b(A)}^A))
\subset T^n.
\end{eqnarray*}
It is a torus of dimension $b(A)$. Finally, observe that the 
$b(A)$-dimensional hypercube 
$\uuuu{p}_j+d^{-1}C(\uuuu{c}_{k_1^A},...,\uuuu{c}_{k_{b(A)}^A})$ 
maps by $\uuuu{e}$ to a subset of $T_A$,
\begin{eqnarray}\label{3.5}
\uuuu{e}\left(\uuuu{p}_j+d^{-1}C(\uuuu{c}_{k_1^A},...,
\uuuu{c}_{k_{b(A)}^A})\right)\subset T_A,
\end{eqnarray}
because the vectors $\uuuu{p}_j$ and $d^{-1}\uuuu{c}_{k_i^A}$
are linear combinations of the vectors $\uuuu{d}_1^A,...,
\uuuu{d}_{b(A)}^A$. And observe that many of the $2^{b(A)}$ vertices
of the image in $T_A$ of this hypercube coincide and that the
vertices form the set $\{\uuuu{e}(\uuuu{p}_j),
\uuuu{e}(\uuuu{p}_{(j+1)_{\mmod d}}),...,
\uuuu{e}(\uuuu{p}_{(j+b(A))_{\mmod d}})\}$, 
because of \eqref{3.4}. 
\end{notations}

Some of these notations are due to Cooper \cite{Co82}, namely 
the deformed simplices $\Delta_A$, the tori $T_A$, 
the {\it block beginnings}, their number $b(A)$ and the 
{\it thick} sets $A$. The hypercubes, the integers 
$\uuuu{a}_k$ and the vectors $\uuuu{b}_k,\uuuu{c}_k,\uuuu{d}_k^A$ and $\uuuu{p}_j$ are new.
The observations in the following lemma are all due to Cooper
\cite{Co82}.

\begin{lemma}\label{t3.2}\cite{Co82}
We stick to the cycle type singularity $f$ above and all induced
data in the Notations \ref{t3.1}.
\begin{eqnarray}\label{3.6}
G&=&\bigcup_{j=1}^d\Delta_N\odot \{\uuuu{e}(\uuuu{p}_j)\}\cup
\bigcup_{A\subsetneqq N,\, A\neq\emptyset}
\Delta_A\odot T_A\\
&=&\bigcup_{j=1}^d\Delta_N\odot \{\uuuu{e}(\uuuu{p}_j)\}\cup
\bigcup_{A\subsetneqq N\textup{ thick}}
\Delta_A\odot T_A.\nonumber
\end{eqnarray}
For $A\subsetneqq N$ with $A\neq\emptyset$, 
\begin{eqnarray}\label{3.7}
G\cap \C^A=\Delta_A\odot T_A.
\end{eqnarray}
The sets $\Delta_N\odot\{\uuuu{e}(\uuuu{p}_j)\}$
and $\Delta_A\odot T_A$ are called {\sf cells} of $G$.
The natural map $\textup{Int}(\Delta_A)\times T_A\to 
\textup{Int}(\Delta_A)\odot T_A$ is a diffeomorphism.
The natural map 
$\Delta_N\to \Delta_N\odot\{\uuuu{e}(\uuuu{p}_j)\}$
is a diffeomorphism. 
\begin{eqnarray}\label{3.8}
&&\dim_\R(\Delta_N\odot\{\uuuu{e}(\uuuu{p}_j)\}=n-1,\\
&&\dim_\R(\Delta_A\odot T_A)=|A|-1+b(A)
\left\{\begin{array}{ll}
=n-1&\textup{if }A\textup{ is thick,}\\
<n-1&\textup{if }A\textup{ is not thick.}\end{array}\right.
\nonumber
\end{eqnarray}
If $A,B\subsetneqq N$ with $A\neq B$, then
$\textup{Int}(\Delta_A)\odot T_A \cap 
\textup{Int}(\Delta_B)\odot T_B=\emptyset$. 
If $B\subsetneqq A\subsetneqq N$ with $b(B)=b(A)$ then
$\Delta_B\odot T_B\subset \Delta_A\odot T_A$. 

$A\subset N$ thick implies $|A|\in\{[\frac{n+1}{2}],
[\frac{n+1}{2}]+1,...,n-1\}$. In the case $n$ even, there are
only two thick sets $A$ with $|A|=\frac{n}{2}$, the set
$A_{od}:=\{1,3,,...,n-1\}$ and the set 
$A_{ev}:=\{2,4,...,n\}$. The cells 
$\Delta_{A_{od}}\odot T_{A_{od}}$ and
$\Delta_{A_{ev}}\odot T_{A_{ev}}$ are the unit spheres
in $\C^{A_{od}}$ respectively $\C^{A_{ev}}$. 

$A\subset N$ almost thick implies
$|A|\in\{[\frac{n}{2}],...,n-2\}$. 
\end{lemma}

The proof is easy. We will not give details. If 
$z_j,z_{j+1}\in\C^*$ with $z_j^{a_j}z_{j+1}>0$ then
$\arg(z_{j+1})\equiv (-a_j)\arg(z_j)\mmod 2\pi.$
This observation is crucial. 

Now we will build up chains, starting with the cells 
$\Delta_N\odot\{\uuuu{e}(\uuuu{p}_j)\}$, and ending
with cycles which represent elements of $H_{n-1}(G,\Z)$.
In section \ref{c4} we will show that these cycles 
(and the cells $\Delta_{A_{od}}\odot T_{A_{od}}$ and
$\Delta_{A_{ev}}\odot T_{A_{ev}}$ in the case $n$ even) 
generate $H_{n-1}(G,\Z)$. 
Of course, it is necessary to refine the cells of $G$
to a simplicial chain complex. We will not describe this
precisely. Except from the cells 
$\Delta_N\odot\{\uuuu{e}(\uuuu{p}_j)\}$, we will work only with
the following chains (compare \eqref{3.5}), 
\begin{eqnarray}\label{3.9}
\CC(A,j):=\Delta_A\odot 
\uuuu{e}\left(\uuuu{p}_j+d^{-1}C(\uuuu{c}_{k_1^A},...,
\uuuu{c}_{k_{b(A)}^A})\right)\subset \Delta_A\odot T_A
\hspace*{0.5cm}\\
\textup{for }A\subsetneqq N\textup{ thick},
\quad j\in\{1,...,d\},\nonumber\\
\CC(B,A,j):=\Delta_B\odot 
\uuuu{e}\left(\uuuu{p}_j+d^{-1}C(\uuuu{c}_{k_1^A},...,
\uuuu{c}_{k_{b(A)}^A})\right)\subset \Delta_B\odot T_A
\hspace*{0.5cm}
\label{3.10}\\
\textup{for }B\subsetneqq A\subsetneqq N\textup{ with }
A\textup{ thick,}\ |B|=|A|-1\geq 1,\ j\in\{1,...,d\}.\nonumber
\end{eqnarray}
Here our notation is not precise in two ways.  (1) If the projection
\begin{eqnarray*}
\uuuu{e}:\uuuu{p}_j+d^{-1}C(\uuuu{c}_{k_1^A},...,
\uuuu{c}_{k_{b(A)}^A})\to T^n
\end{eqnarray*}
is not injective, the chain $\CC(A,j)$ or $\CC(B,A,j)$
shall take multiplicities into account. 
(2) The chain obtains an orientation from the order of the
vectors $\uuuu{c}_{k_1^A},...,\uuuu{c}_{k_{b(A)}^A}$
if they are linearly independent. 

We decompose the boundary of a chain $\CC(A,j)$ 
into two parts, $\paa\CC(A,j)=\paa_1\CC(A,j)+\paa_2\CC(A,j)$
with
\begin{eqnarray}\label{3.11}
\paa_1\CC(A,j)&=& \paa\Delta_A\odot 
\uuuu{e}\left(\uuuu{p}_j+d^{-1}C(\uuuu{c}_{k_1^A},...,
\uuuu{c}_{k_{b(A)}^A})\right),\\
\paa_2\CC(A,j)&=& (-1)^{|A|-1}\Delta_A\odot 
\paa \uuuu{e}\left(\uuuu{p}_j+d^{-1} C(\uuuu{c}_{k_1^A},...,
\uuuu{c}_{k_{b(A)}^A})\right)\nonumber\\
&=& (-1)^{|A|-1}\Delta_A\odot 
\uuuu{e}\left(\uuuu{p}_j+d^{-1}\paa C(\uuuu{c}_{k_1^A},...,
\uuuu{c}_{k_{b(A)}^A})\right).\hspace*{0.5cm}\label{3.12}
\end{eqnarray}

The definitions of the following chains $R_j^{(k)}$ and
$X_j^{(k)}$ are crucial. Here $j\in\{1,...,d\}$
and $k\in\{0,1,...,[\frac{n+1}{2}]\}$. 
\begin{eqnarray}\label{3.13}
R_j^{(k)}&:=& \sum_{A\textup{ thick},|A|=n-k}\sign(A)\cdot\CC(A,j)
\quad\textup{for }k\in\{1,2,...,[\frac{n}{2}]\},\hspace*{1cm}\\
R_j^{(\frac{n+1}{2})}&:=&0\quad\textup{ if }n\textup{ is odd,}
\nonumber\\
R_j^{(0)}&:=& \Delta_N\odot \{\uuuu{e}(\uuuu{p}_j)\}.\nonumber
\nonumber
\end{eqnarray}
For $R_j^{(k)}$ with $k\geq 1$, the decomposition $\paa R_j^{(k)}
=\paa_1 R_j^{(k)}+\paa_2 R_j^{(k)}$ of its boundary 
into two parts is well defined. 
For $k=0$ we define $\paa_2 R_j^{(0)}:=0$ and $\paa_1 R_j^{(0)}:=
\paa R_j^{(0)}$. 
The definition of the next chains is inductive.
Again $j\in\{1,...,d\}$.
\begin{eqnarray}\label{3.14}
X_j^{(0)}&:=& \Delta_N\odot \{\uuuu{e}(\uuuu{p}_j)\}=R_j^{(0)},\\
X_j^{(k)}&:=& X_j^{(k-1)}-X_{j+1}^{(k-1)}+R_j^{(k)}
\quad\textup{for }k\in\{1,...,[\frac{n+1}{2}]\}.\label{3.15}
\end{eqnarray}

\begin{theorem}\label{t3.3}
Again, $j\in\{1,...,d\}$ and $k\in\{0,1,...,[\frac{n+1}{2}]\}$.
\begin{eqnarray}\label{3.16}
\paa X_j^{(k)}&=& \paa_1 R_j^{(k)},\\ 
\textup{and especially }\quad\paa X_j^{([\frac{n+1}{2}])}&=&0.
\label{3.17}
\end{eqnarray}
So, the chains $X_j^{([\frac{n+1}{2}])}$ are cycles and
induce homology classes in $H_{n-1}(G,\Z)$. 
(Theorem \ref{t4.6} will tell more.)
If $n$ is odd then
\begin{eqnarray}\label{3.18}
\sum_{j=1}^d X_j^{(\frac{n+1}{2})}=0.
\end{eqnarray}
If $n$ is even then
\begin{eqnarray}\label{3.19}
\sum_{j=1}^d X_j^{(\frac{n}{2})}&=&
(-1)^{(\frac{n}{2}+2)(\frac{n}{2}+1)\frac{1}{2}}\cdot 
\uuuu{a}_1\uuuu{a}_3...\uuuu{a}_{n-1}\nonumber \\
&\cdot& \left(\Delta_{A_{od}}\odot T_{A_{od}} + (-1)^{\frac{n}{2}}
a_1a_3...a_{n-1}\Delta_{A_{ev}}\odot T_{A_{ev}}\right).
\hspace*{1cm}
\end{eqnarray}
\end{theorem}

{\bf Proof:} \eqref{3.16} for $k=0$ is 
$\paa X_j^{(0)}=\paa_1 R_j^{(0)}$. It is true by definition
of $\paa_1 R_j^{(0)}$. By induction, we obtain for $k\geq 0$
\begin{eqnarray*}
\paa X_j^{(k+1)}&=& \paa X_j^{(k)}-\paa X_{j+1}^{(k)}
+\paa R_j^{(k+1)}\\
&=& \left(\paa_1 R_j^{(k)}-\paa_1 R_{j+1}^{(k)}+\paa_2 R_j^{(k+1)}
\right) + \paa_1 R_j^{(k+1)}.
\end{eqnarray*}
Therefore we have to show for $k\geq 0$ 
\begin{eqnarray}\label{3.20}
\paa_1 R_j^{(k)}-\paa_1 R_{j+1}^{(k)}+\paa_2 R_j^{(k+1)}=0.
\end{eqnarray}

In fact, $R_j^{(k+1)}$ was chosen so that \eqref{3.20} holds.
In order to show \eqref{3.20}, first we study 
$\paa_1 R_j^{(k)}$. 

For $A=\{k_1,...,k_{|A|}\}\subset N$ with $k_1<...<k_{|A|}$
if $n\notin A$ and with $k_1=n\ \&\ k_2<...<k_{|A|}$ if $n\in A$
and for $B=\{k_1,...,\widehat{k_j},...,k_{|A|}\}\subset A$ define
\begin{eqnarray}\label{3.21}
\sign(B,A):=(-1)^{j-1}.
\end{eqnarray}
Then
\begin{eqnarray}
\paa_1 R_j^{(k)}&=& \sum_{A\textup{ thick},|A|=n-k}
\sum_{B\subset A,|B|=|A|-1}\sign(B,A)\cdot\sign(A)\cdot\CC(B,A,j)
\nonumber \\
&&\textup{for }k\geq 1,\label{3.22}\\
\paa_1 R_j^{(0)}&=& \sum_{B:\, |B|=n-1}\sign(B,N)\cdot
\Delta_B\odot \{\uuuu{e}(\uuuu{p}_j)\}.\nonumber
\end{eqnarray}

First we consider the cases $k\geq 1$. A part of the following
arguments will be valid also for the case $k=0$ and will
give this case. In general, the subsets $B$ of thick sets $A$
with $|B|=|A|-1=n-k-1$ are of three different types,
\begin{eqnarray*}
\begin{array}{l|l|l}
\textup{type I} & \textup{type II} & \textup{type III}\\
\hline b(B)=b(A)+1 & b(B)=b(A) & b(B)=b(A)-1\\
\Rightarrow\, B\textup{ thick}&
\Rightarrow\, B\textup{ almost thick}&
\end{array}
\end{eqnarray*}

First, consider a set $B$ of type III.
Then only one thick set $A$ with $|A|=n-k$ and
$A\supset B$ exists. 
It contains a block $\{k_j^A\}$ which consists of
a single number $k_j^A$, and $B=A-\{k_j^A\}$. 
The gap of $B$ which contains $k_j^A$ consists of $k_j^A$
and the gaps of $A$ left and right of $k_j^A$, so this gap
of the set $B$ has 3 elements. Therefore the set
$B$ is neither thick nor almost thick. 
All other gaps of $B$ consist of a single number. Observe
\begin{eqnarray*}
\pr_B(\uuuu{c}_{k_j^A})&=&\pr_B(\uuuu{c}_{k_{j+1}^A}),\quad
\textup{thus}\\\ 
\pr_B(d^{-1}C(\uuuu{c}_{k_1^A},...,\uuuu{c}_{k_{b(A)}^A}))
&=& \pr_B(d^{-1}C(\uuuu{c}_{k_1^B},...,\uuuu{c}_{k_{b(B)}^B})),
\quad\textup{thus}\\
\CC(B,A,j)&=& \CC(B,j)\subset\Delta_B\odot T_B,\\
\dim\CC(B,A,j)&=& |B|-1+b(B)
=|A|-2+b(A)-1\\ &=&\dim\CC(A,j)-2=\dim R_j^{(k)}-2.
\end{eqnarray*}
Therefore this part $\CC(B,A,j)$ of the boundary 
$\paa_1 R_j^{(k)}$ has too small dimension and can be ignored. 

Next, consider a set $B$ of type II.
It is almost thick. 
It has one gap which consists of two numbers $k_0$ and
$(k_0+1)_{\mmod n}$. All other gaps consist of a single number.
The only two thick sets $A$ with $A\supset B$ and $|A|=n-k$ are 
$B^{(1)}:=B\cup\{k_0\}$ and $B^{(2)}:=B\cup\{(k_0+1)_{\mmod n}\}$.
All possibilities for $k_0$ except one are easy to treat,
namely the cases $k_0\in\{1,...,n\}-\{n-1\}$. 
In all these cases
\begin{eqnarray*}
\pr_B(\uuuu{c}_{k_i^{B^{(1)}}})=\pr_B(\uuuu{c}_{k_i^{B^{(2)}}})
=\pr_B(\uuuu{c}_{k_i^{B}})\textup{ for any }i\in\{1,...,b(B)\},\\
\textup{thus}\quad 
\CC(B,B^{(1)},j)=\CC(B,B^{(2)},j)=\CC(B,j).
\end{eqnarray*}
One checks also
\begin{eqnarray}\label{3.23}
\sign(B,B^{(1)})=\sign(B,B^{(2)}), \quad 
\sign(B^{(1)})=-\sign(B^{(2)}).
\end{eqnarray}
These observations show that these contributions to 
$\paa_1 R_j^{(k)}$ cancel.

The only difficult possibility for $k_0$ is the case $k_0=n-1$. 
Then
\begin{eqnarray*}
\begin{array}{ccccccccc}
k_1^{B^{(1)}}=1&<&k_2^{B^{(1)}}&<&...&<&k_{b(B)}^{B^{(1)}} ,
&, & k_j^{B^{(1)}}=k_j^B\\
 & & ll & & & & ll & & \\
 & & k_1^{B^{(2)}} &<&...&<&k_{b(B)-1}^{B^{(2)}}
 &<&k_{b(B)}^{B^{(2)}}=n ,
\end{array}
\end{eqnarray*}
\begin{eqnarray}
\uuuu{c}_{k_{b(B)}^{B^{(2)}}}=
\uuuu{c}_n&=&(-1)^nd\cdot \uuuu{b}_n+\uuuu{c}_1,\nonumber\\
\pr_B(\uuuu{c}_{k_{b(B)}^{B^{(2)}}})&=&
((-1)^nd+1)\pr_B(\uuuu{c}_1),\nonumber \\
1=\sign(B,B^{(2)})&=& (-1)^{|B|}\sign(B,B^{(1)}),\label{3.24}\\
\sign(B^{(2)})&=& (-1)^{n-1}\sign(B^{(1)}),\nonumber \\
\sign(B,B^{(2)})\sign(B^{(2)})&=& (-1)^{b(B)}
\sign(B,B^{(1)})\sign(B^{(1)}).\nonumber
\end{eqnarray}
This together with \eqref{3.3}, formulas for $\uuuu{c}_k$ 
and $\uuuu{d}_j^B$ and especially $k_1^B=1$ and 
$\uuuu{b}_1=(0,...,0)$ shows 
\begin{eqnarray*}
&& \pr_T\pr_B(\uuuu{p}_j+d^{-1}C(\uuuu{c}_{k_1^{B^{(2)}}},...,
\uuuu{c}_{k_{b(B)}^{B^{(2)}}})) \\
&&+(-1)^{b(B)}\pr_T\pr_B
(\uuuu{p}_j+d^{-1}C(\uuuu{c}_{k_1^{B^{(1)}}},...,
\uuuu{c}_{k_{b(B)}^{B^{(1)}}}))\\
&=& \pr_T\pr_B(\uuuu{p}_j+d^{-1}C(\uuuu{c}_{k_2^{B}},
...,\uuuu{c}_{k_{b(B)}^B},(-1)^nd \uuuu{c}_1))\\
&=& \pr_Tpr_B(\uuuu{p}_j+C(d^{-1}\uuuu{c}_{k_2^{B}},
...,d^{-1}\uuuu{c}_{k_{b(B)}^B},(-1)^n\uuuu{c}_1))\\
&\stackrel{\eqref{3.3}}{=}& \pr_T\pr_B(\uuuu{p}_j+C( 
d^{-1}\uuuu{c}_{k_2^{B}}-d^{-1}\uuuu{c}_1,...,
d^{-1}\uuuu{c}_{k_{b(B)}^B}-d^{-1}\uuuu{c}_1,(-1)^n\uuuu{c}_1))\\
&=& \pr_T\pr_B(\uuuu{p}_j+C(
(-1)^n\uuuu{b}_{k_2^{B}},...,
(-1)^n\uuuu{b}_{k_{b(B)}^B},(-1)^n\uuuu{c}_1))\\
&=& (-1)^{nb(B)} 
\uuuu{a}_1\uuuu{a}_{k_2^B}...\uuuu{a}_{k_{b(B)}^B}\cdot
\pr_T\pr_B
(C(\uuuu{d}_1^B,\uuuu{d}_2^B,...,\uuuu{d}_{b(B)}^B)).
\end{eqnarray*}
This implies 
\begin{eqnarray}
&&\CC(B,B^{(2)},j)
+(-1)^{b(B)}\cdot \CC(B,B^{(1)},j)\nonumber\\
&=& (-1)^{nb(B)}
\uuuu{a}_1\uuuu{a}_{k_2^B}...\uuuu{a}_{k_{b(B)}^B}\cdot
\Delta_B\odot T_B.\label{3.25}
\end{eqnarray}
This is up to the sign $\sign(B,B^{(2)})\cdot\sign(B^{(2)})$
the contribution of $B$ to $\paa_1 R_j^{(k)}$. 
Especially, it is independent of $j$.
Therefore a set $B$ of type II makes the same contribution
to $\paa_1 R_j^{(k)}$ and to $\paa_1 R_{j+1}^{(k)}$. 
Thus its contribution to the difference
$\paa_1 R_j^{(k)}-\paa_1 R_{j+1}^{(k)}$ is zero.

Finally, consider a set $B$ of type I. It is thick.
Its gaps are the sets $\{(k_1^B-1)_{\mmod n}\},...,
\{(k_{b(B)}^B-1)_{\mmod n}\}$. Exactly $b(B)$ thick
sets $A$ with $A\supset B$ and $|A|=n-k$ exist.
They are the sets $A^{(i)}=B\cup\{(k_i^B-1)_{\mmod n}\}$
for $i\in\{1,...,b(B)\}$. 

The set of block beginnings of $A^{(i)}$ is the set
$\{k_1^B,...,\widehat{k_i^B},...,k_{b(B)}^B\}$. Therefore
the contribution of $B$ to the boundary $\paa_1 R_j^{(k)}$ is
\begin{eqnarray}
\sum_{i=1}^{b(B)}&&\sign(B,A^{(i)})\cdot\sign(A^{(i)})\nonumber\\
&\cdot&\Delta_B\odot \uuuu{e}\left(
\uuuu{p}_j+d^{-1}C(\uuuu{c}_{k_1^B},...,
\widehat{\uuuu{c}_{k_i^B}},...,\uuuu{c}_{k_{b(B)}^B})\right),
\label{3.26}
\end{eqnarray}
and the contribution of $B$ to $\paa_1 R_{j+1}^{(k)}$ looks
analogously, with $j$ replaced by $j+1$. 
The following calculation of signs will be useful,
\begin{eqnarray}
\sign(B)&=& (-1)^{(i-1)+(n-k_i^B)+(b(B)-i)}\cdot \sign(A^{(i)})
\nonumber\\
&=&(-1)^{k_i^B+|B|-1}\cdot\sign(A^{(i)}),\nonumber \\
\sign(B,A^{(i)})&=&(-1)^{(k_i^B-1)-(i-1)}=(-1)^{k_i^B-i},
\nonumber\\
\sign(B,A^{(i)})\cdot\sign(A^{(i)})
&=& \sign(B)\cdot(-1)^{|B|+i-1}.
\label{3.27}
\end{eqnarray}

On the other side, the boundary $\paa_2 R_j^{(k+1)}$ has only
contributions from sets $B$ of type I. The contribution
of one such set $B$ is as follows.
Here we use \eqref{3.12}, \eqref{3.2} and \eqref{3.4}.
\begin{eqnarray*}
&&\textup{contribution of }B\textup{ to }\paa_2 R_j^{(k+1)}\\
&\stackrel{\eqref{3.12}}{=}&\sign(B)(-1)^{|B|-1}\cdot 
\Delta_B\odot \uuuu{e}\left(\uuuu{p}_j+d^{-1}
\paa C(\uuuu{c}_{k_1^B},...,\uuuu{c}_{k_{b(B)}^B})\right)\\
&\stackrel{\eqref{3.2}}{=}& \sign(B)(-1)^{|B|-1}\cdot \\
&&\Bigl(
\sum_{i=1}^{b(B)}(-1)^{i-1}\Delta_B\odot\uuuu{e}\left(\uuuu{p}_j + 
d^{-1}C(\uuuu{c}_{k_1^B},...,
\widehat{\uuuu{c}_{k_i^B}},...,\uuuu{c}_{k_{b(B)}^B})\right) 
\Bigr.\\
&&\Bigl. 
-\sum_{i=1}^{b(B)}(-1)^{i-1}\Delta_B\odot\uuuu{e}\left(\uuuu{p}_j 
+ d^{-1}\uuuu{c}_{k_i^B} + 
d^{-1}C(\uuuu{c}_{k_1^B},...,
\widehat{\uuuu{c}_{k_i^B}},...,\uuuu{c}_{k_{b(B)}^B})\right)
\Bigr)
\end{eqnarray*}
\begin{eqnarray*}
&\stackrel{\eqref{3.4}}{=}& \sum_{i=1}^{b(B)}\sign(B)(-1)^{|B|+i}
\Delta_B\odot\uuuu{e}\left(\uuuu{p}_j + 
d^{-1}C(\uuuu{c}_{k_1^B},...,
\widehat{\uuuu{c}_{k_i^B}},...,\uuuu{c}_{k_{b(B)}^B})\right)\\
&&- \sum_{i=1}^{b(B)}\sign(B)(-1)^{|B|+i}
\Delta_B\odot\uuuu{e}\left(\uuuu{p}_{j+1}
+d^{-1}C(\uuuu{c}_{k_1^B},...,
\widehat{\uuuu{c}_{k_i^B}},...,\uuuu{c}_{k_{b(B)}^B})\right)\\
&=& -(\textup{contribution of }B\textup{ to }\paa_1 R_j^{(k)})\\
&& +(\textup{contribution of }B\textup{ to }\paa_1 R_{j+1}^{(k)}).
\end{eqnarray*}
Therefore the contribution of $B$ to 
$\paa_1 R_j^{(k)}-\paa_1 R_{j+1}^{(k)} + \paa_2 R_j^{(k+1)}$
is zero. 

So, no set $B\subset N$ with $|B|=n-k-1$ gives a contribution
to this sum. Therefore this sum is zero. \eqref{3.20} is proved
for $k\geq 1$. 

Now we will show \eqref{3.20} for $k=0$.  
$\paa_1 R_j^{(0)}$, $\paa_1 R_{j+1}^{(0)}$ and $\paa_2 R_j^{(1)}$ are sums over sets $B\subset N$ with $|B|=n-1$. 
These sets are the sets $B^{(i)}:=N-\{i\}$ for $i\in N$. 
The set $B^{(i)}$ is thick with $b(B^{(i)})=1$. 
It can be treated as the sets of type I above. 
The same calculations as above for sets of type I show
that the contribution of $B^{(i)}$ to the 
sum $\paa_1 R_j^{(0)}-\paa_1 R_{j+1}^{(0)} +\paa_2 R_j^{(1)}$ 
is zero. Therefore this sum is zero. Therefore \eqref{3.20} holds
also for $k=0$, so it holds for all $k\geq 0$. 

This implies \eqref{3.16} for $k\geq 1$. It holds for $k=0$
anyway. 

It remains to prove \eqref{3.17}, \eqref{3.18} and \eqref{3.19}.
\eqref{3.17} is because of \eqref{3.16} equivalent to
$\paa_1 R_j^{([\frac{n+1}{2}])}=0$. 
In the case $n$ odd this is trivial as then
$R_j^{(\frac{n+1}{2})}=0$ by definition. In the case $n$ even,
we have to show $\paa_1 R_j^{(\frac{n}{2})}=0$. 
By Lemma \ref{t3.2}, there are only two thick sets $A$ with
$|A|=\frac{n}{2}$, the sets $A_{od}=\{1,3,...,n-1\}$
and $A_{ev}=\{2,4,...,n\}$. In the discussion above of
$\paa_1 R_j^{(k)}$ for $k\geq 1$, we distinguished 3 different
types of sets $B$ with $|B|=|A|-1=n-k-1$ and $B\subset A$.
In the cases $A\in \{A_{od},A_{ev}\}$, we have only sets $B$ 
of type III. Above we stated that the part $\CC(B,A,j)$ of 
$\paa_1 R_j^{(k)}$ (here for $k=\frac{n}{2}$) has too 
small dimension and can be ignored. Therefore
$\paa_1 R_j^{(\frac{n}{2})}=0$. We proved \eqref{3.17}.

\eqref{3.18} holds for odd $n$ and is an immediate consequence
of \eqref{3.15} and $R_j^{(\frac{n+1}{2})}=0$.

\eqref{3.19} holds for even $n$. It requires two calculations 
which are similar to the one which led to \eqref{3.25}. 
The details are as follows. \eqref{3.15} and \eqref{3.13} show
\begin{eqnarray*}
\sum_{j=1}^d X_j^{(\frac{n}{2})}
&=& \sum_{j=1}^d R_j^{(\frac{n}{2})}\\
&=& \sum_{j=1}^d\sign(A_{od})\cdot\CC(A_{od},j)
+ \sum_{j=1}^d\sign(A_{ev})\cdot\CC(A_{ev},j)\\
&=& \sum_{j=1}^d \sign(A_{od})\cdot
\Delta_{A_{od}}\odot \uuuu{e}\left(\uuuu{p}_j
+d^{-1}C(\uuuu{c}_1,\uuuu{c}_3,...,\uuuu{c}_{n-1})\right)\\
&+& \sum_{j=1}^d \sign(A_{ev})\cdot
\Delta_{A_{ev}}\odot \uuuu{e}\left(\uuuu{p}_j
+d^{-1}C(\uuuu{c}_2,\uuuu{c}_4,...,\uuuu{c}_{n})\right).
\end{eqnarray*}
Here
\begin{eqnarray*}
\sign(A_{od})&=&(-1)^{\frac{n}{2}(\frac{n}{2}+1)\frac{1}{2}},
\nonumber\\ 
\sign(A_{ev})&=&(-1)^{(\frac{n}{2}-1)\frac{n}{2}\frac{1}{2}}
=(-1)^{\frac{n}{2}}\cdot \sign(A_{od}).
\end{eqnarray*}
Because of $\uuuu{p}_{j+1}=\uuuu{p}_j+d^{-1}\uuuu{c}_1$
and \eqref{3.3} (and $(-1)^n=1$ as $n$ is even), we have
\begin{eqnarray*}
&& \sum_{j=1}^d \pr_T\left(
\uuuu{p}_j+d^{-1} C(\uuuu{c}_1,\uuuu{c}_3,...,\uuuu{c}_{n-1})
\right)\\
&=& \pr_T\left(\uuuu{p}_1+C(\uuuu{c}_1,d^{-1}\uuuu{c}_3,....,
d^{-1}\uuuu{c}_{n-1})\right)\\
&\stackrel{\eqref{3.3}}{=}& 
\pr_T\left(\uuuu{p}_1+C(\uuuu{c}_1,d^{-1}\uuuu{c}_3
-d^{-1}\uuuu{c}_1,...,d^{-1}\uuuu{c}_{n-1}-d^{-1}\uuuu{c}_1)
\right)\\
&=& \pr_T\left(\uuuu{p}_1+C(\uuuu{c}_1,(-1)^n\uuuu{b}_3,...,
(-1)^n\uuuu{b}_{n-1})\right)\\
&=& (-1)^{n/2-1}\pr_T\left(C(\uuuu{b}_3,...,\uuuu{b}_{n-1},
\uuuu{c}_1)\right).
\end{eqnarray*}
If we identify $\pr_{A_{od}}(\uuuu{b}_j)$ with a column vector
in $M_{n/2\times 1}(\Z)$ then the tuple 
$(\pr_{A_{od}}(\uuuu{b}_3),...,\pr_{A_{od}}(\uuuu{b}_{n-1}),
\pr_{A_{od}}(\uuuu{c}_1))$ is an upper triangular 
$\frac{n}{2}\times\frac{n}{2}$-matrix, and its determinant is
$\uuuu{a}_1\uuuu{a}_3...\uuuu{a}_{n-1}.$
Therefore
\begin{eqnarray*}
&&\sum_{j=1}^d\Delta_{A_{od}}\odot\uuuu{e}\left(
\uuuu{p}_j+d^{-1}C(\uuuu{c}_1,\uuuu{c}_3,...,\uuuu{c}_{n-1})
\right)\\
&=& (-1)^{n/2-1}\uuuu{a}_1\uuuu{a}_3...\uuuu{a}_{n-1}\cdot
\Delta_{A_{od}}\odot T_{A_{od}}.
\end{eqnarray*}
Observe
$$\pr_{A_{ev}}(\uuuu{c}_{2l})=\pr_{A_{ev}}(\uuuu{c}_{2l-1}).$$
Therefore 
\begin{eqnarray*}
&&\sum_{j=1}^d\pr_{A_{ev}}\pr_T\left(\uuuu{p}_j
+d^{-1}C(\uuuu{c}_2,\uuuu{c}_4,...,\uuuu{c}_{n})\right)\\
&=& \sum_{j=1}^d\pr_{A_{ev}}\pr_T\left(\uuuu{p}_j
+d^{-1}C(\uuuu{c}_1,\uuuu{c}_3,...,\uuuu{c}_{n-1})\right)\\
&=& (-1)^{n/2-1}\pr_{A_{ev}}\pr_T\left(
C(\uuuu{b}_3,...,\uuuu{b}_{n-1},\uuuu{c}_1)\right).
\end{eqnarray*}
If we identify $\pr_{A_{ev}}(\uuuu{b}_j)$ with a column vector
in $M_{n/2\times 1}(\Z)$ then the tuple 
$(\pr_{A_{ev}}(\uuuu{b}_3),...,\pr_{A_{ev}}(\uuuu{b}_{n-1}),
\pr_{A_{ev}}(\uuuu{c}_1))$ is an upper triangular 
$\frac{n}{2}\times\frac{n}{2}$-matrix, and its determinant is
$$\uuuu{a}_2\uuuu{a}_4...\uuuu{a}_{n}=
\uuuu{a}_1\uuuu{a}_3...\uuuu{a}_{n-1}\cdot a_1a_3...a_{n-1}.$$
Therefore
\begin{eqnarray*}
&&\sum_{j=1}^d\Delta_{A_{ev}}\odot\uuuu{e}\left(
\uuuu{p}_j+d^{-1}C(\uuuu{c}_2,\uuuu{c}_4,...,\uuuu{c}_{n})
\right)\\
&=& (-1)^{n/2-1}\uuuu{a}_1\uuuu{a}_3...\uuuu{a}_{n-1}\cdot
a_1a_3...a_{n-1}\cdot \Delta_{A_{ev}}\odot T_{A_{ev}}.
\end{eqnarray*}
\eqref{3.19} follows. \hfill$\Box$

\section{A spectral sequence, following and correcting Cooper}\label{c4}
\setcounter{equation}{0}

\noindent
Throughout this section we fix 
a cycle type singularity $f(x_1,...,x_n)$ as in
\eqref{1.2} with $n\geq 2$ and $a_1,...,a_n\in\N$ with
\eqref{1.3}, as in section \ref{c3}. 
Cooper \cite{Co82} considered a filtration of the set 
$G\subset \oooo{F_0}=f^{-1}(\R_{\geq 0})\cap S^{2n-1}
\subset S^{2n-1}\subset\C^n$, which was defined in \eqref{3.1}, 
by a sequence of subsets,
$$G=G_n\supset G_{n-1}\supset ...\supset G_1\supset G_0=
\emptyset.$$
He studied a spectral sequence 
$(E^r_{s,t},d^r_{s,t})_{r\geq 1,s\geq 1,s+t\geq 0}$ which is
associated to this filtration. It allows to determine 
$H_{n-1}(G,\Z)$.
Though he made two serious mistakes. Here we will partly 
follow his line of thoughts, but correct the mistakes.
Remark \ref{t5.1} will explain the differences to \cite{Co82}.

Cooper does not give a reference for the construction of the
spectral sequence. We found a manuscript of Hatcher \cite{Ha04}
very useful and will cite from it. 

First we will define the subspaces $G_s\subset G$
and make some elementary observations. 
Then we will introduce the spectral sequence, following
Hatcher. Then we will study it in detail (following
partly Cooper) in Lemma \ref{t4.3} and Theorem \ref{t4.5}.
Theorem \ref{t4.6} will give the conclusion for $H_{n-1}(G,\Z)$. 
It complements Theorem \ref{t3.3}.

\begin{definition}\label{t4.1}\cite{Co82}
Recall the definition \eqref{3.1} of the set $G\subset \C^n$.
Define $G_n:=G$, $G_0:=\emptyset$, and define for
$s\in\{1,...,n-1\}$
\begin{eqnarray}
G_s&:=&\{z\in G\,|\, \textup{at most }s
\textup{ of the coordinates } z_1,...,z_n\textup{ are }\neq 0\}
\nonumber \\
&=& \bigcup_{A\subset N:\, |A|\leq s}\Delta_A\odot T_A.\label{4.1}
\end{eqnarray}
\end{definition}

Here the second equality follows from Lemma \ref{t3.2}.
Lemma \ref{t3.2} gives also the differences $G_s-G_{s-1}$,
\begin{eqnarray}\label{4.2}
G_n-G_{n-1}&=& \bigcup_{j=1}^d\textup{Int}(\Delta_N)\odot 
\{\uuuu{e}(\uuuu{p}_j)\},\\
G_s-G_{s-1}&=&  \bigcup_{A\subset N:\, |A|=s}
\textup{Int}(\Delta_A)\odot T_A\textup{ for }
s\in\{1,...,n-1\}.\label{4.3}
\end{eqnarray}

In this paper, all considered homology groups will have
coefficients in $\Z$. 
Especially, we consider for $(s,t)\in N\times\Z$ 
the homology groups of the spaces $G_s$ and of the pairs 
$(G_s,G_{s-1})$,
\begin{eqnarray}\label{4.4}
A^1_{s,t}&:=& H_{s+t}(G_s):=H_{s+t}(G_s,\Z),\\
E^1_{s,t}&:=& H_{s+t}(G_s,G_{s-1}):=H_{s+t}(G_s,G_{s-1},\Z).
\label{4.5}
\end{eqnarray}
We extend this definition by 
\begin{eqnarray*}\nonumber
A^1_{s,t}=E^1_{s,t}=0&&\textup{for }s\leq 0,\\
A^1_{s,t}=H_{s+t}(G),\ E^1_{s,t}=0&&\textup{for }
s\geq n+1.
\end{eqnarray*}
Of course
\begin{eqnarray}\label{4.6}
A^1_{s,t}=E^1_{s,t}=0&&\textup{for }s+t<0\textup{ or }s+t\geq n,
\end{eqnarray}
as $\dim_\R G_s\leq n-1$.
Cooper \cite{Co82} observed that Lemma \ref{t3.2}
(especially that the map 
$\textup{Int}(\Delta_A)\times T_A\to
\textup{Int}(\Delta_A)\odot T_A$ is a diffeomorphism)
and \eqref{4.3} imply
\begin{eqnarray}\label{4.7}
H_*(G_n,G_{n-1})&=& \bigoplus_{j=1}^d 
H_{n-1}(\Delta_N\odot\{\uuuu{e}(\uuuu{p}_j)\},
\paa\Delta_N\odot\{\uuuu{e}(\uuuu{p}_j)\}),\\
H_*(G_s,G_{s-1})&=& \bigoplus_{A\subset N:\, |A|=s}
H_*(\Delta_A\odot T_A,\paa\Delta_A\odot T_A)\label{4.8}\\
&&\hspace*{3cm} \textup{ for }s\leq n-1.\nonumber
\end{eqnarray}
And here
\begin{eqnarray}
H_{n-1}(\Delta_N\odot\{\uuuu{e}(\uuuu{p}_j)\},
\paa\Delta_N\odot\{\uuuu{e}(\uuuu{p}_j)\})\cong\Z\nonumber\\
\textup{ with generator the class }[X_j^{(0)}],\label{4.9}\\
H_*(\Delta_A\odot T_A,\paa\Delta_A\odot T_A)
\cong H_{|A|-1}(\Delta_A,\paa \Delta_A)\otimes 
H_*(T_A),\label{4.10}\\
H_{|A|-1}(\Delta_A,\paa\Delta_A)\cong\Z
\textup{ with generator the class }[\Delta_A],\label{4.11}\\
H_*(T_A)\cong \bigotimes_{j=1}^{b(A)} H_*(S^1),\quad
H_*(S^1)=H_0(S^1)\oplus H_1(S^1),\label{4.12}\\
H_0(S^1)\cong\Z\ \textup{ and }\ H_1(S^1)\cong\Z\nonumber\\
\textup{ with generators }[p]\textup{ and }[S^1]
\textup{ for any }p\in S^1.\label{4.13}
\end{eqnarray}
By \eqref{4.7} and \eqref{4.9} the group
\begin{eqnarray}\label{4.14}
E^1_{n,-1}=H_{n-1}(G_n,G_{n-1})=H_*(G_n,G_{n-1})
=\bigoplus_{j=1}^d \Z\cdot [X_j^{(0)}]
\end{eqnarray}
is a $\Z$-lattice of rank $d$ with generators the classes
$[X_1^{(0)}],...,[X_d^{(0)}]$.
All the groups $E^1_{s,t}$ are $\Z$-lattices
(= finitely generated free $\Z$-modules). 

The long exact homology sequence for the pair
$(G_s,G_{s-1})$ for $s\in N$ reads as follows,
\begin{eqnarray*}
\cdots A^1_{s-1,t+1}\stackrel{i^1}{\to} A^1_{s,t}
\stackrel{j^1}{\to} E^1_{s,t}
\stackrel{k^1}{\to}A^1_{s-1,t}
\stackrel{i^1}{\to}A^1_{s,t-1}
\stackrel{j^1}{\to}E^1_{s,t-1}
\stackrel{k^1}{\to}A^1_{s-1,t-1} \cdots
\end{eqnarray*}
Here $i^1$ is induced by the embedding 
$G_{s-1}\hookrightarrow G_s$, $j^1$ is the natural
map from absolute homology to relative homology, 
and $k^1$ is the boundary map.
Together, these exact homology sequences can be put into
a large diagram which Hatcher \cite{Ha04} calls a 
{\it staircase diagram}, because the long exact sequences look
like staircases in this diagram:
\begin{eqnarray*}
\begin{array}{ccccccccc}
A^1_{s-1,t+1}& \stackrel{j^1}{\to} & E^1_{s-1,t+1} &
\stackrel{k^1}{\to} & A^1_{s-2,t+1} & \stackrel{j^1}{\to} &
E^1_{s-2,t+1} & \stackrel{k^1}{\to} & A^1_{s-3,t+1} \\
\downarrow i^1 & & & & \downarrow i^1 & & & & \downarrow i^1 \\
A^1_{s,t}& \stackrel{j^1}{\to} & E^1_{s,t} &
\stackrel{k^1}{\to} & A^1_{s-1,t} & \stackrel{j^1}{\to} &
E^1_{s-1,t} & \stackrel{k^1}{\to} & A^1_{s-2,t} \\
\downarrow i^1 & & & & \downarrow i^1 & & & & \downarrow i^1 \\
A^1_{s+1,t-1}& \stackrel{j^1}{\to} & E^1_{s+1,t-1} &
\stackrel{k^1}{\to} & A^1_{s,t-1} & \stackrel{j^1}{\to} &
E^1_{s,t-1} & \stackrel{k^1}{\to} & A^1_{s-1,t-1}
\end{array}
\end{eqnarray*}
The map
\begin{eqnarray}\label{4.15}
d^1_{s,t}:=d^1:=j^1\circ k^1:E^1_{s,t}\to E^1_{s-1,t}
\end{eqnarray} 
is the cellular boundary map \cite{Ha04}.
The system $(E^1_{s,t},d^1_{s,t})_{s,t\in\Z}$ of spaces
and maps is the first page of a spectral sequence which
is associated to the filtration $G=G_n\supset G_{n-1}\supset
...\supset G_1\supset G_0=\emptyset$. 
The other pages $(E^r_{s,t},d^r_{s,t})_{s,t\in\Z}$ for 
$r\geq 2$ are constructed inductively together with 
versions of the staircase diagram above for any $r\geq 2$
instead of $r=1$, so, tuples 
$(A^r_{s,t},E^r_{s,t},i^r,j^r,k^r,d^r)$ are constructed 
inductively for any $r\geq 2$. The following theorem
describes this construction and gives the general properties.
It follows from Lemma 5.1 and Proposition 5.2 in
\cite{Ha04}. Proposition 5.2 in \cite{Ha04} applies with
the conditions (i) and (ii) in it. 
The indices in $A^1_{s,t}$ and $E^1_{s,t}$ here are chosen
differently (more standard) from those in \cite{Ha04}.

\begin{theorem}\label{t4.2} 
\cite[Lemma 5.1 and Proposition 5.2]{Ha04}

(a) (Properties) Fix $r\geq 2$. The tuple   
$(A^r_{s,t},E^r_{s,t},i^r,j^r,k^r,d^r)$ will be constructed
from the tuple 
$(A^{r-1}_{s,t},E^{r-1}_{s,t},i^{r-1},j^{r-1},k^{r-1},d^{r-1})$
in part (c). It has the following properties:
\begin{eqnarray}\label{4.16}
i^r:A^r_{s,t}&\to & A^r_{s+1,t-1},\\
j^r:A^r_{s,t}&\to & E^r_{s-r+1,t+r-1},\label{4.17}\\
k^r:E^r_{s,t}&\to & A^r_{s-1,t},\label{4.18}\\
d^r:=j^r\circ k^r:E^r_{s,t}&\to&  E^r_{s-r,t+r-1},
\quad d^r\circ d^r=0.\label{4.19}
\end{eqnarray}
The restriction of $i^r$ to $A^r_{s,t}$ is also called
$i^r_{s,t}$, and similarly for $j^r,k^r$ and $d^r$. 
There are long exact sequences,
\begin{eqnarray*}
\cdots A^r_{s+r-2,t-r+2}\stackrel{i^r}{\to} A^r_{s+r-1,t-r+1}
\stackrel{j^r}{\to} E^r_{s,t}
\stackrel{k^r}{\to}A^r_{s-1,t}\\
\stackrel{i^r}{\to}A^r_{s,t-1}
\stackrel{j^r}{\to}E^r_{s-r+1,t+r-2}
\stackrel{k^r}{\to}A^r_{s-r,t+r-2} \cdots
\end{eqnarray*}
They can be put together into the following staircase diagram
(the positions of the arrows for $j^r$ are not precise).
\begin{eqnarray*}
\begin{array}{ccccccccc}
A^r_{s-1,t+1}& \stackrel{\nearrow}{j^r} & E^r_{s-1,t+1} &
\stackrel{k^r}{\to} & A^r_{s-2,t+1} & \stackrel{\nearrow}{j^r} &
E^r_{s-2,t+1} & \stackrel{k^r}{\to} & A^r_{s-3,t+1} \\
\downarrow i^r & & & & \downarrow i^r & & & & \downarrow i^r \\
A^r_{s,t}& \stackrel{\nearrow}{j^r} & E^r_{s,t} &
\stackrel{k^r}{\to} & A^r_{s-1,t} & \stackrel{\nearrow}{j^r} &
E^r_{s-1,t} & \stackrel{k^r}{\to} & A^r_{s-2,t} \\
\downarrow i^r & & & & \downarrow i^r & & & & \downarrow i^r \\
A^r_{s+1,t-1}& \stackrel{\nearrow}{j^r} & E^r_{s+1,t-1} &
\stackrel{k^r}{\to} & A^r_{s,t-1} & \stackrel{\nearrow}{j^r} &
E^r_{s,t-1} & \stackrel{k^r}{\to} & A^r_{s-1,t-1}
\end{array}
\end{eqnarray*}

\medskip
(b) (Spectral sequence)
The part $(E^r_{s,t},d^r_{s,t})_{r\geq 1,s,t\in\Z}$ of the 
tuples above for all $r\geq 1$ is a spectral sequence. 
It converges to $H_*(G)$. This imprecise statement means the
following. For each $(s,t)\in\Z^2$, a bound $r(s,t)\in\N$ exists
such that the space $E^\infty_{s,t}:=E^{r(s,t)}_{s,t}$
coincides with $E^r_{s,t}$ for any $r\geq r(s,t)$, and
\begin{eqnarray}\label{4.20} 
E^\infty_{s,t}&\cong & F^s_{s+t}/F^{s-1}_{s+t},\textup{ where}\\
F^s_{s+t}&:=& \Imm(A^1_{s,t}\to A^1_{n,s+t-n})
=\Imm(H_{s+t}(G_s)\to H_{s+t}(G))\nonumber.
\end{eqnarray}

\medskip
(c) (Construction) Fix $r\geq 2$. 
Suppose that the tuple 
$(A^{r-1}_{s,t},E^{r-1}_{s,t},i^{r-1},j^{r-1},k^{r-1},d^{r-1})$
has been constructed. Then especially $d^{r-1}\circ d^{r-1}=0$.
Therefore the quotient
\begin{eqnarray}\label{4.21}
E^r_{s,t}&:=& \ker(d^{r-1}_{s,t})/\Imm(d^{r-1}_{s+r-1,t-r+2})
\end{eqnarray}
is well-defined. It is a subquotient of $E^{r-1}_{s,t}$.
The space
\begin{eqnarray}\label{4.22}
A^r_{s,t}&:=& i^{r-1}(A^{r-1}_{s-1,t+1})\subset A^{r-1}_{s,t}
\end{eqnarray}
is well-defined, anyway. It is a subspace of $A^{r-1}_{s,t}$. 
The map $i^r$ is the restriction of $i^{r-1}$ to $A^r_{s,t}$.
The map $j^r$ on $A^r_{s,t}=i^{r-1}(A^{r-1}_{s-1,t+1})$ is
defined on $i^{r-1}(a)$ for $a\in A^{r-1}_{s-1,t+1}$ by
$j^r(i^{r-1}(a)):=[j^{r-1}(a)]$. A priori, this
is in $E^{r-1}_{s-r+1,t+r-1}/\Imm(d^{r-1}_{s,t+1})$. 
But it turns out to be in $E^r_{s-r+1,t+r-1}$.
The map $k^r$ on $E^r_{s,t}$ is defined on $[a]$ for 
$a\in \ker(d^{r-1}_{s,t})\subset E^{r-1}_{s,t}$ by
$k^r([a]):=k^{r-1}(a)$ which is a priori in 
$A^{r-1}_{s-1,t}$. But it
turns out to be in $A^r_{s-1,t}$ and to be well-defined. 
\end{theorem}

Theorem \ref{t4.2} gives basic and general properties,
which are not specific to the geometry of $G$ and its
subspaces $G_s$ (except for our definition of 
$F^s_{s+t}$ where we used $A^1_{n,s+t-n}=H_{s+t}(G)$).
Now we will use the geometry of $G$ and its subspaces
$G_s$ to make more specific statements.
We are interested especially in $E^\infty_{s,n-1-s}$
as that is the subquotient $F^s_{n-1}/F^{s-1}_{n-1}$ 
of $H_{n-1}(G)$ by part (b). 
Lemma \ref{t4.3} gives first elementary observations.
Algebraic statements in Lemma \ref{t4.4} make a part 
of Theorem \ref{t4.5} more transparent. 
Theorem \ref{4.5} states the main results on the 
spectral sequence.
Theorem \ref{t4.6} gives the conclusion for
$H_{n-1}(G)$. It complements Theorem \ref{t3.3}.

\begin{lemma}\label{t4.3}
(a) For $r\geq 2$ 
\begin{eqnarray}\label{4.23}
E^r_{s,n-1-s}=\ker(d^{r-1}_{s,n-1-s})\subset E^{r-1}_{s,n-1-s}
\subset ...\subset E^1_{s,n-1-s}.
\end{eqnarray}
(b) For $r\geq 1$ and for $s<[\frac{n+1}{2}]$
\begin{eqnarray}\label{4.24}
E^r_{s,n-1-s}=E^\infty_{s,n-1-s}=0\quad
\textup{and thus }F^s_{n-1}=0.
\end{eqnarray}
(c) Consider even $n$. For $r\geq 1$
\begin{eqnarray}\nonumber
F^{n/2}_{n-1}&=&\Imm\left((H_{n-1}(G_{\frac{n}{2}})\to H_{n-1}(G)
\right) \cong E^\infty_{\frac{n}{2},\frac{n}{2}-1}\\
&=&E^r_{\frac{n}{2},\frac{n}{2}-1}=E^1_{\frac{n}{2},\frac{n}{2}-1}
= H_{n-1}(G_{\frac{n}{2}},G_{\frac{n}{2}-1})\cong\Z^2.\label{4.25}
\end{eqnarray}
The generators are the classes in $H_{n-1}(G)$ respectively
in $H_{n-1}(G_{\frac{n}{2}},G_{\frac{n}{2}-1})$ of the spheres
$\Delta_{A_{od}}\odot T_{A_{od}}\subset \C^{A_{od}}$ and 
$\Delta_{A_{ev}}\odot T_{A_{ev}}\subset \C^{A_{ev}}$
(see Lemma \ref{t3.2}). 
\end{lemma}

{\bf Proof:}
(a) Because of \eqref{4.6} $E^{r-1}_{s,t}=0$ for $s+t\geq n$. 
Therefore $d^{r-1}_{s+r-1,n-1-s-r+2}=0$, and by \eqref{4.21}
$E^r_{s,n-1-s}=\ker(d^{r-1}_{s,n-1-s})$.

(b) For $s\leq n-1$, the set $E^1_{s,n-1-s}=H_{n-1}(G_s,G_{s-1})$
has by \eqref{4.8}, \eqref{4.10} and \eqref{4.12}
only contributions from sets $A\subset N$ with
$b(A)=n-|A|$, so from thick sets $A$. By Lemma \ref{t3.2},
if $s<[\frac{n+1}{2}]$, there are no thick sets $A$ with
$|A|=s$. Therefore then $E^1_{s,n-1-s}=0$
and $E^r_{s,n-1-s}=0$ for $r\geq 1$ and $E^\infty_{s,n-1-s}=0$.
This and \eqref{4.20} imply inductively $F^s_{n-1}=0$ 
for $s<[\frac{n+1}{2}]$.

(c) The first equality is the definition of 
$F^{n/2}_{n-1}$. The isomorphism $F^{n/2}_{n-1}
\cong E^\infty_{\frac{n}{2},\frac{n}{2}-1}$ follows with 
\eqref{4.20} and $F^s_{n-1}=0$ for $s<\frac{n}{2}$.

The equality
$E^1_{\frac{n}{2},\frac{n}{2}-1}
= H_{n-1}(G_{\frac{n}{2}},G_{\frac{n}{2}-1})$ is the definition
of $E^1_{\frac{n}{2},\frac{n}{2}-1}$. By \eqref{4.8},
\eqref{4.10} and \eqref{4.12} it has a contribution isomorphic
to $\Z$ with generator $[\Delta_A\odot T_A]$ for any thick
set $A$ with $|A|=\frac{n}{2}$. By Lemma \ref{t3.2},
these are only the sets $A_{od}$ and $A_{ev}$. 
Therefore $E^1_{\frac{n}{2},\frac{n}{2}-1}\cong\Z^2$.
As the spheres $\Delta_{A_{od}}\odot T_{A_{od}}$ and 
$\Delta_{A_{ev}}\odot T_{A_{ev}}$ have no boundary,
the boundary map 
\begin{eqnarray*}
k^1_{\frac{n}{2},\frac{n}{2}-1}: E^1_{\frac{n}{2},\frac{n}{2}-1}
\to A^1_{\frac{n}{2}-1,\frac{n}{2}-1} 
=H_{n-2}(G_{\frac{n}{2}-1})
\end{eqnarray*}
and also the induced maps $k^r_{\frac{n}{2},\frac{n}{2}-1}$ 
and the maps  $d^r_{\frac{n}{2},\frac{n}{2}-1}$ 
are zero maps. Therefore for $r\geq 1$ 
\begin{eqnarray*}
E^1_{\frac{n}{2},\frac{n}{2}-1}
=E^r_{\frac{n}{2},\frac{n}{2}-1}
=E^\infty_{\frac{n}{2},\frac{n}{2}-1}.
\end{eqnarray*}
Part (c) is proved. \hfill$\Box$

\bigskip
The sole purpose of the following elementary algebraic lemma
is to make Theorem \ref{t4.5} (d) more transparent. 

\begin{lemma}\label{t4.4}
Recall from \eqref{4.14} that the group 
$E^1_{n,-1}=H_{n-1}(G,G_{n-1})$ is a $\Z$-lattice of rank $d$
with generators the classes $[X_1^{(0)}],...,[X_d^{(0)}]$.
Define for $r\in\{1,...,[\frac{n+3}{2}]\}$ the sublattice
\begin{eqnarray}\label{4.26}
\www{E}^r_{n,-1}:=\sum_{j=1}^d\Z\cdot [X_j^{(r-1)}]
\subset E^1_{n,-1},
\end{eqnarray}
which is generated by the classes 
$[X_1^{(r-1)}],...,[X_d^{(r-1)}]$ in $H_{n-1}(G,G_{n-1})$. 
Then $\www{E}^1_{n,-1}=E^1_{n,-1}$. For $r\geq 2$ 
\begin{eqnarray}\label{4.27}
[X_j^{(r-1)}]&=&[X_j^{(r-2)}]-[X_{(j+1)_{\mmod d}}^{(r-2)}]
\quad\textup{for }j\in\{1,...,d\},\\
\sum_{j=1}^d[X_j^{(r-1)}]&=&0,\label{4.28}
\end{eqnarray}
and $\www{E}^r_{n,-1}\subset E^1_{n,-1}$ is a sublattice
of rank $d-1$ which is generated by any $d-1$ of the $d$
elements $[X_1^{(r-1)}],...,[X_d^{(r-1)}]$.
The lattice $\www{E}^2_{n,-1}$ is  a primitive sublattice 
of $E^1_{n,-1}$.
For $r\geq 3$, $\www{E}^r_{n,-1}$ has index $d$ in 
$\www{E}^{r-1}_{n,-1}$. 
\end{lemma}

{\bf Proof:}
\eqref{4.27} follows from the definition of $X_j^{(r-1)}$ 
in \eqref{3.15} and from $R^{(r-1)}_j\subset G_{n-1}$.
\eqref{4.28} is an immediate consequence of \eqref{4.27}. 
$\www{E}^2_{n,-1}$ is obviously a primitive sublattice of
$E^1_{n,-1}$ of rank $d-1$. For $r\geq 3$, 
$\www{E}^r_{n,-1}$ has index $d$ in $\www{E}^{r-1}_{n,-1}$
because
\begin{eqnarray}\label{4.29}
\www{E}^{r}_{n,-1} = \{\sum_{j=1}^dz_k\cdot [X_j^{(r-2)}]\,|\, 
z_k\in\Z, \sum_{j=1}^d z_j\in d\Z\}
\end{eqnarray}
by \eqref{4.27} and \eqref{4.28} (for $[X_j^{(r-2)}]$).
\hfill $\Box$

\begin{theorem}\label{t4.5}
(a) \cite[18.]{Co82} 
For $s\in\{[\frac{n+2}{2}],[\frac{n+4}{2}],...,n-1\}$
\begin{eqnarray}\label{4.30}
d^1_{s,n-1-s}:E^1_{s,n-1-s}\to E^1_{s-1,n-1-s}
\quad\textup{is injective,}\\
E^r_{s,n-1-s}=E^\infty_{s,n-1-s}=0\qquad\textup{for }r\geq 2.\label{4.31}
\end{eqnarray}
(b) For $t\in\{3,...,[\frac{n+3}{2}]\}$ and $3\leq r\leq t$
\begin{eqnarray}\label{4.32}
E^{r-1}_{n-t+1,t-3}\subset 
\frac{E^1_{n-t+1,t-3}}{d^1(E^1_{n-t+2,t-3})}.
\end{eqnarray}
(c) For $s<[\frac{n+2}{2}]$ and $r\geq 1$
\begin{eqnarray}\label{4.33}
E^r_{s-1,n-1-s}=0.
\end{eqnarray}
(d) Recall the definition of $\www{E}^r_{n,-1}$ in Lemma 
\ref{t4.4}.
\begin{eqnarray}\label{4.34}
E^r_{n,-1}&=& \www{E}^r_{n,-1} \quad\textup{for }
r\in\{1,...,[\frac{n+3}{2}]\},\\
E^\infty_{n,-1}&=&E^{r}_{n,-1}=\www{E}^{[\frac{n+3}{2}]}_{n,-1}
\quad\textup{for }r\geq [\frac{n+3}{2}].\label{4.35}
\end{eqnarray}
\end{theorem}

{\bf Proof:}
(a) \eqref{4.31} is an immediate consequence of \eqref{4.30},
because of the definition of $E^r_{s,n-1-s}$ in \eqref{4.21}.

The injectivity of $d^1_{s,n-1-s}$ was proved in
\cite[18.]{Co82}. Our proof differs in a way which allows
to apply it also to the proof of \eqref{4.34} in part (d). 
First we consider some useful sets and maps.
Fix $s\in\{[\frac{n+2}{2}],[\frac{n+4}{2}],...,n-1\}$.
\begin{eqnarray*}
\AAA_1(s)&:=& \{A\subset N\,|\, A\textup{ thick},|A|=s\},\\
\AAA_2(s)&:=& \{(A,j)\,|\, A\in\AAA_1(s), j\in\{1,...,b(A)\}\\
&&\hspace*{3.5cm}\textup{ with }k_j^A,(k_j^A+1)_{\mmod n}\in A\},
\\
\pr_1&:&\AAA_2(s)\to\AAA_1(s),(A,j)\mapsto A,
\textup{ the canonical projection}.
\end{eqnarray*}
For the given $s$, each $A\in\AAA_1(s)$ contains a block which
consists of $\geq 2$ elements. Therefore 
$\pr_1:\AAA_2(s)\to \AAA_1(s)$ is surjective. 
\begin{eqnarray*}
\BB(s)&:=& \{B\subset N\,|\, B\textup{ almost thick},
|B|=s-1\},\\
\beta_1&:&\BB(s)\to N,\beta_1(B):=k_0\textup{ if }
\{k_0,(k_0+1)_{\mmod n}\}\subset N-B,\\
\beta_2&:&\BB(s)\to N,\beta_2(B):=(\beta_1(B)+1)_{\mmod n},\\
\alpha_1&:&\AAA_2(s)\to\BB(s),(A,j)\mapsto 
A-\{(k_{(j+1)_{\mmod d}}-2)_{\mmod n}\},\\
\alpha_2&:&\AAA_2(s)\to\BB(s),(A,j)\mapsto 
A-\{k_j^A\}.
\end{eqnarray*}
$\beta_1(B)$ and $\beta_2(B)$ are the first and last element
of the unique block of $N-B$ with two elements.
$k_j^A$ is the beginning of a block of $A$
with $\geq 2$ elements, and $(k_{(j+1)_{\mmod d}}-2)_{\mmod n}\}$
is the last element of this block.
$\alpha_1$ and $\alpha_2$ are bijections, and
(see the notations \ref{t3.1} for $B^{(1)}$ and $B^{(2)}$)
\begin{eqnarray*}
B^{(1)}=B\cup\{\beta_1(B)\}=\pr_1(\alpha_1^{-1}(B)),\\
B^{(2)}=B\cup\{\beta_2(B)\}=\pr_1(\alpha_2^{-1}(B)).
\end{eqnarray*}

Recall what \eqref{4.8}, \eqref{4.10} and \eqref{4.12}
say about $E^1_{s,n-1-s}$ and $E^1_{s-1,n-1-s}$.
Both are $\Z$-lattices. The generators of 
$E^1_{s,n-1-s}$ are simply the classes $[\Delta_A\odot T_A]$,
$A\in\AAA_1(s)$, 
\begin{eqnarray}\label{4.36}
&&E^1_{s,n-1-s}= H_{n-1}(G_s,G_{s-1})\\
&&\hspace*{1cm} =
\bigoplus_{A\in \AAA_1(s)}H_{n-1}(\Delta_A\odot T_A,
\paa\Delta_A\odot T_A)
=\bigoplus_{A\in \AAA_1(s)}\Z\cdot [\Delta_A\odot T_A],
\nonumber
\end{eqnarray}
$E^1_{s-1,n-1-s}$ splits into two parts, one from 
almost thick sets, the other from thick sets,
\begin{eqnarray}\label{4.37}
&& E^1_{s-1,n-1-s}=H_{n-2}(G_{s-1},G_{s-2})\\
&&\hspace*{1cm}= \bigl(H_{n-2}(G_{s-1},G_{s-2})\bigr)_{\BB}
\oplus \bigl(H_{n-2}(G_{s-1},G_{s-2})\bigr)_{\AAA}
\textup{ with }\nonumber\\
&& \bigl(H_{n-2}(G_{s-1},G_{s-2})\bigr)_{\BB} 
:= \bigoplus_{B\in \BB(s)}\Z\cdot [\Delta_B\odot T_B],
\nonumber\\
&&\bigl(H_{n-2}(G_{s-1},G_{s-2})\bigr)_{\AAA}
:= \bigoplus_{A\in \AAA_1(s-1)}
H_{s-1}(\Delta_A,\paa\Delta_A)\otimes H_{b(A)-1}(T_A)\nonumber\\
&&\hspace*{4cm}=\bigoplus_{A\in \AAA_1(s-1)}\bigoplus_{i=1}^{b(A)}\Z\cdot 
[\Delta_A\odot T_{A\cup\{(k_i^A-1)_{\mmod n}\}}].\nonumber
\end{eqnarray}
We will work mainly with the part from almost thick sets. 
The projection to this part is called $\pr_\BB$,
\begin{eqnarray*}
\pr_\BB:E^1_{s-1,n-1-s}=H_{n-2}(G_{s-1},G_{s-2})
\to \left(H_{n-2}(G_{s-1},G_{s-2})\right)_\BB.
\end{eqnarray*}

We will prove that $\pr_\BB\circ d^1_{s,n-1-s}$ is injective.
This implies that $d^1_{s,n-1-s}$ is injective. 
Recall $d^1_{s,n-1-s}=j^1_{s-1,n-1-s}\circ k^1_{s,n-1-s}$,
and $k^1_{s,n-1-s}$ is a boundary map. 
Therefore for $A\in\AAA_1(s)$ 
\begin{eqnarray}\label{4.38}
&&d^1_{s,n-1-s}([\Delta_A\odot T_A])
= [\paa\Delta_A\odot T_A] \\
&&\hspace*{1cm}= \sum_{B\subset A:\, B\in\BB(s)\cup\AAA_1(s-1)}
\sign(B,A)\cdot [\Delta_B\odot T_A],\nonumber
\end{eqnarray}
We care only about the terms for $B\subset A$ with 
$B\in\BB(s)$. These sets $B$ split into two types,
\begin{eqnarray}\label{4.39}
\{B\subset A\,|\, B\in\BB(s)\} 
=\{B\subset A\,|\, A=B^{(1)}\}
\cup \{B\subset A\,|\, A=B^{(2)}\}.
\end{eqnarray}
We claim to have on the level of chains
\begin{eqnarray}\label{4.40}
\Delta_B\odot T_A =\left\{\begin{array}{ll}
\Delta_B\odot T_B&\textup{if }A=B^{(1)},\\
(-a_{\beta_2(B)})\cdot \Delta_B\odot T_B & 
\textup{if }A=B^{(2)},\beta_2(B)\neq n,\\
(-1)^{b(B)-1}(-a_n) \Delta_B\odot T_B & 
\textup{if }A=B^{(2)},\beta_2(B)=n.\\
\end{array}\right.
\end{eqnarray}
This follows from comparison of $T_A$ with $T_B$. 
If $B=\alpha_1(A,j)$ then 
$\beta_1(B)=(k_{(j+1)_{\mmod n}}-2)_{\mmod n}$
is the last element of the $j$-th block of $A$.
If $B=\alpha_2(A,j)$ then $\beta_2(B)=k_j^A$
is the beginning of the $j$-th block of $A$. 
The generating vectors of $T_A$ and $T_B$ from these
blocks of $A$ and $B$ differ as follows, 
\begin{eqnarray}
\pr_B(\uuuu{d}_j^A) =
\pr_B(\uuuu{d}_j^B)&\textup{if }B=\alpha_1(A,j),\nonumber\\
\pr_B(\uuuu{d}_j^A) =
(-a_{\beta_2(B)})\cdot \pr_B(\uuuu{d}_{\www{j}}^B) & 
\textup{if }B=\alpha_2(A,j),\label{4.41}
\end{eqnarray}
where $\www{j}=j$ if $\beta_2(B)\neq n$ and where 
in the case $\beta_2(B)=n$ 
\begin{eqnarray*}
j=b(B)=b(A),\ \www{j}=1,\\
\pr_B(\uuuu{d}_i^A)=\pr_B(\uuuu{d}_{i+1}^B)\textup{ for }
i\in\{1,...,b(B)-1\}.
\end{eqnarray*}
This shows \eqref{4.40}. 
Together \eqref{4.38}, \eqref{4.39} and \eqref{4.40}
give
\begin{eqnarray}
&&\pr_\BB\left(d^1_{s,n-1-s}([\Delta_A\odot T_A])\right)
= \sum_{B\in \alpha_1(\pr_1^{-1}(A))}
\sign(B,A)\cdot [\Delta_B\odot T_B] \nonumber\\
&&+ \sum_{B\in \alpha_2(\pr_1^{-1}(A)),\beta_2(B)\neq n}
\sign(B,A)\cdot (-a_{\beta_2(B)})\cdot [\Delta_B\odot T_B]
\label{4.42}\\
&&+ \sum_{B\in \alpha_2(\pr_1^{-1}(A)),\beta_2(B)=n}
(-1)^{b(B)-1}
\sign(B,A)\cdot (-a_n)\cdot [\Delta_B\odot T_B].
\nonumber
\end{eqnarray}
Recall from \eqref{3.23} and \eqref{3.24}  
\begin{eqnarray*}
\sign(B,B^{(2)})&=&\sign(B,B^{(1)})\quad
\textup{ if }\beta_2(B)\neq n,\\
1=\sign(B,B^{(2)})&=&(-1)^{|B|}\sign(B,B^{(1)})\quad\textup{ if }
\beta_2(B)=n,\\
\textup{and }(-1)^{b(B)-1+|B|}&=&(-1)^n.
\end{eqnarray*}
In view of these signs and \eqref{4.42},  
an arbitrary linear combination 
$\sum_{A\in\AAA_1(s)}z_A\cdot [\Delta_A\odot T_A]$,
$z_A\in\Z$, is mapped by $\pr_\BB\circ d^1_{s,n-1-s}$ to
\begin{eqnarray}\label{4.43}
&&\pr_\BB\bigl(d^1_{s,n-1-s}\left(
\sum_{A\in\AAA_1(s)}z_A\cdot [\Delta_A\odot T_A]\right)\bigr)\\
&=&\sum_{B\in \BB(s),\beta_2(B)\neq n}
\sign(B,B^{(1)})\bigl(z_{B^{(1)}}-a_{\beta_2(B)}
z_{B^{(2)}}\bigr)[\Delta_B\odot T_B]\Bigr. \nonumber\\
&+&\Bigl.\sum_{B\in \BB(s),\beta_2(B)=n}
\sign(B,B^{(1)})\bigl(z_{B^{(1)}}-(-1)^{n}a_n
z_{B^{(2)}}\bigr)[\Delta_B\odot T_B].
\nonumber
\end{eqnarray}
The following Claim 1 will be useful here and in the proof
of \eqref{4.34} in part (d).

\medskip
{\bf Claim 1:} {\it Fix any $B\in \BB(s)$. A sequence
$(B_i)_{i\in N}$ of elements $B_i\in\BB(s)$ with $B=B_n$ and 
with the following property exists,
\begin{eqnarray}\label{4.44}
B^{(1)}_{(i+1)_{\mmod n}} = B^{(2)}_i\quad\textup{for }i\in N.
\end{eqnarray}
Additionally, 
\begin{eqnarray}\label{4.45}
\textup{the map}\quad 
N\to N,\ i\mapsto \beta_2(B_i),\quad\textup{is a bijection.}
\end{eqnarray}
}

\medskip
{\bf Proof of Claim 1:} 
We will describe such a sequence from a point of view
which will make its existence clear. Suppose we have such
a sequence $(B_i)_{i\in N}$. Compare the set of gaps of
$B_{(i-1)_{\mmod n}}^{(2)}=B_{i}^{(1)}$ with the set of gaps of 
$B_{i}^{(2)}=B_{(i+1)_{\mmod n}}^{(1)}$.
These sets almost coincide. Only the gap 
$\{\beta_2(B_{i})\}$ of $B_i^{(1)}$ is shifted one position
to the left to the gap $\{\beta_1(B_i)\}$ of $B_i^{(2)}$. 
Starting from the set of gaps of $B_n^{(2)}=B_1^{(1)}$,
one shifts in $n$ steps all gaps to the left, so that the
final position of each gap is the original position of the
gap left of it. One has to take care that the gaps stay always
apart from one another. 
It is clear that starting from the set of gaps
of $B_n^{(2)}$, one can find such $n$ steps. Also 
\eqref{4.45} is clear. This finishes the proof 
of Claim 1. \hfill$(\Box)$

\medskip
Consider now a linear combination 
$\sum_{A\in\AAA_1(s)}z_A\cdot [\Delta_A\odot T_A]$,
$z_A\in\Z$, which is mapped by $\pr_\BB\circ d^1_{s,n-1-s}$ to
$0$. Choose any $A\in\AAA_1(s)$ and choose $B\in\BB(s)$
with $B^{(2)}=A$. Claim 1 provides a sequence
$(B_i)_{i\in N}$ with $B_n=B$, $B_n^{(2)}=A$, 
\eqref{4.44} and \eqref{4.45}.
Then \eqref{4.43} gives the relations
\begin{eqnarray}\label{4.46}
z_{B_{(i-1)_{\mmod n}}^{(2)}} = z_{B_i^{(1)}}
=\left\{\begin{array}{ll}
a_{\beta_2(B_i)}z_{B_i^{(2)}}\quad\textup{if }
\beta_2(B_i)\neq n,\\
(-1)^{n}a_nz_{B_i^{(2)}}\quad\textup{if }
\beta_2(B_i)=n.\end{array}\right. 
\end{eqnarray}
Together they give because of \eqref{4.45} especially
\begin{eqnarray*}
z_A=\uuuu{a}_{n+1}\cdot z_A,
\quad \textup{so }(-1)^nd\cdot z_A=0,
\quad \textup{so } z_A=0.
\end{eqnarray*}
Therefore $\pr_\BB\circ d^1_{s,n-1-s}$ is injective,
and thus also $d^1_{s,n-1-s}$ is injective. 
This finishes the proof of part (a).

(b) The definition \eqref{4.21} of $E^{r-1}_{n-t+1,t-3}$ for
$r\geq 3$ gives 
\begin{eqnarray*}
E^{r-1}_{n-t+1,t-3}&\subset&
\frac{E^{r-2}_{n-t+1,t-3}}{d^{r-2}(E^{r-2}_{n-t+r-1,t-r})}.
\end{eqnarray*}
If $r=3$, this is \eqref{4.32}. If $r\geq 4$, \eqref{4.31}
gives $E^{r-2}_{n-t+r-1,t-r}=0$, so 
$E^{r-1}_{n-t+1,t-3}\subset E^{r-2}_{n-t+1,t-3}$.
Induction gives \eqref{4.32}.

(c) $E^1_{s-1,n-1-s}=H_{n-2}(G_{s-1},G_{s-2})$ has by
\eqref{4.8}, \eqref{4.10} and \eqref{4.12} only contributions
from sets $A\subset N$ with $|A|=s-1$ and $b(A)\in\{n-|A|,
n-1-|A|\}$. These are thick or almost thick sets.
By Lemma \ref{t3.2}, for $s<[\frac{n+2}{2}]$, there are no
thick or almost thick sets with $|A|=s-1$. Therefore then
$E^1_{s-1,n-1-s}=0$ and $E^r_{s-1,n-1-s}=0$ for $r\geq 1$. 

(d) By \eqref{4.23}, for $r\geq 2$ 
\begin{eqnarray}\label{4.47}
E^r_{n,-1}=\ker(d^{r-1}_{n,-1}:E^{r-1}_{n,-1}\to 
E^{r-1}_{n-r+1,r-3})\subset E^{r-1}_{n,-1}\subset E^1_{n,-1}.
\end{eqnarray}
For $r\geq[\frac{n+3}{2}]+1$, $E^{r-1}_{n-r+1,r-3}=0$
by part (c). Therefore then $E^r_{n,-1}=E^{r-1}_{n,-1}$.
Inductively we obtain
\begin{eqnarray*}
E^\infty_{n,-1}=E^r_{n,-1}=E^{[\frac{n+3}{2}]}_{n,-1}
\quad\textup{for }r\geq [\frac{n+3}{2}].
\end{eqnarray*}

It remains to prove \eqref{4.34}. We will prove it by
induction in $r$. By definition $\www{E}^1_{n,-1}=E^1_{n,-1}$.
First we treat the special case $r=2$:
\begin{eqnarray*}
E^2_{n,-1}&=&\ker(d^1:E^1_{n,-1}\to E^1_{n-1,-1}),\\
d^1([X_j^{(0)}])&=&d^1([\Delta_N\odot\{\uuuu{e}(\uuuu{p}_j)\}])
= [\paa \Delta_N\odot\{\uuuu{e}(\uuuu{p}_j)\}]\\
&=&\sum_{B\subset N:\, |B|=n-1}\sign(B,N)\cdot 
[\Delta_B\odot\{\uuuu{e}(\uuuu{p}_j)\}]\\
&=&\sum_{B\subset N:\, |B|=n-1}\sign(B,N)\cdot 
[\Delta_B\odot\{\uuuu{e}(\uuuu{p}_1)\}]\\
&\in& E^1_{n-1,-1}=H_{n-2}(G_{n-1},G_{n-2})\\
&=&\bigoplus_{i=1}^n
\Z\cdot [\Delta_{N-\{i\}}\odot
\{\uuuu{e}(\uuuu{p}_1)\}]\cong\Z^n.
\end{eqnarray*}
The last equality follows from \eqref{4.8} and \eqref{4.10}.
Obviously $d^1([X_j^{(0)}])$
is independent of $j$ and is $\neq 0$. Therefore 
\begin{eqnarray*}
E^2_{n,-1}&=&\ker(d^1:E^1_{n,-1}\to E^1_{n-1,-1})\\
&=&\{\sum_{j=1}^dz_j\cdot [X_j^{(0)}]\,|\, z_j\in\Z,
\sum_{j=1}^dz_j=0\}=\www{E}^2_{n,-1}.
\end{eqnarray*}

Now we suppose $r\in\{3,...,[\frac{n+3}{2}]\}$ and
(induction hypothesis) $E^{r-1}_{n,-1}=\www{E}^{r-1}_{n,-1}$.
This induction hypothesis and \eqref{4.26} give 
\begin{eqnarray}\label{4.48}
E^{r-1}_{n,-1}=\www{E}^{r-1}_{n,-1}
=\bigoplus_{j=1}^{d-1}\Z\cdot[X_j^{(r-2)}]
=\sum_{j=1}^d\Z\cdot[X_j^{(r-2)}]\subset E^1_{n,-1}.
\end{eqnarray}
We have to control $d^{r-1}_{n,-1}([X_j^{(r-2)}])$. 
Recall four points:
\begin{list}{}{}
\item[(i)] 
$d^{r-1}_{n,-1}=j^{r-1}_{n-1,-1}\circ k^{r-1}_{n,-1}$.
\item[(ii)]
$k^{r-1}_{n,-1}$ is a boundary map with image in
\begin{eqnarray*}
A^{r-1}_{n-1,-1}=(i^1)^{r-2}(A^1_{n-r+1,r-3})
=\Imm(H_{n-2}(G_{n-r+1})\to H_{n-2}(G_{n-1})).
\end{eqnarray*}
\item[(iii)]
$j^{r-1}_{n-1,-1}$ maps this space to
$E^{r-1}_{n-r+1,r-3}$, which satisfies because of
\eqref{4.32} (for $r=t$)
\begin{eqnarray*}
E^{r-1}_{n-r+1,r-3}\subset 
\frac{E^1_{n-r+1,r-3}}{d^1(E^1_{n-r+2,r-3})}
=\frac{H_{n-2}(G_{n-r+1},G_{n-r})}{d^1(E^1_{n-r+2,r-3})}.
\end{eqnarray*}
\item[(iv)]
\eqref{3.16} in Theorem \ref{t3.3} gives 
$\paa X_j^{(r-2)}=\paa_1R_j^{(r-2)}$, and this is a chain in
$G_{n-r+1}$, which fits. 
\end{list}

These four points imply
\begin{eqnarray}\label{4.49}
d^{r-1}_{n,-1}([X_j^{(r-2)}] = [\paa_1 R_j^{(r-2)}]
\in\frac{H_{n-2}(G_{n-r+1},G_{n-r})}{d^1(E^1_{n-r+2,r-3})}.
\end{eqnarray} 
We will show the following claim. 

\medskip
{\bf Claim 2}: {\it The class $[\paa_1R_j^{(r-2)}]$
in the quotient in \eqref{4.49} of the chain $\paa_1R_j^{(r-2)}$
is independent of $j$, and for $m\in\Z$ 
\begin{eqnarray}\label{4.50}
m\cdot [\paa_1R_j^{(r-2)}]=0\iff d|m.
\end{eqnarray}}

\medskip
Claim 2 implies
\begin{eqnarray*}
E^r_{n,-1}&=&\ker(d^{r-1}_{n,-1}:E^{r-1}_{n,-1}
\to E^{r-1}_{n-r+1,r-3})\\
&=&\{\sum_{j=1}^dz_j\cdot [X_j^{(r-2)}]\,|\, z_j\in\Z,
\sum_{j=1}^dz_j\in d\Z\}
\stackrel{\eqref{4.29}}{=}\www{E}^r_{n,-1}.
\end{eqnarray*}
It rests to prove Claim 2. 

The chain $\paa_1R_j^{(r-2)}$ was already studied in the 
proof of Theorem \ref{t3.3}. Formula \eqref{3.22} showed
contributions from pairs $(B,A)$ of sets $B\subset A\subset N$
with $|B|=|A|-1$, $|A|=n-r+2$ and $A$ thick. A priori
$B$ could by of one of three types I, II and III.
But a set $B$ of type III gives nothing.  A set $B$ of type II 
is almost thick and gives either nothing or the contribution 
in \eqref{3.25} (up to a sign) which is independent of $j$. 
A set $B$ of type I is thick and gives the contribution in 
\eqref{3.26}, which we call now 
$\textup{Cont}(B,\paa_1R_j^{(r-2)})$. 

First we will show that
$\textup{Cont}(B,\paa_1R_j^{(r-2)})$ 
is a relative cycle, so that it gives a class 
\begin{eqnarray}\label{4.51}
[\textup{Cont}(B,\paa_1R_j^{(r-2)})]
\in H_{n-2}(G_{n-r+1},G_{n-r}).
\end{eqnarray}
Then we will show that this class is independent of $j$. 

$\paa_1\left(\textup{Cont}(B,\paa_1R_j^{(r-2)})\right)$ contains
$\paa\Delta_B$, so it is in $G_{n-r}$, and we can ignore it.
For the signs in \eqref{3.26}, recall \eqref{3.27}. Then
\begin{eqnarray*}
&&\paa_2\left(\textup{Cont}(B,\paa_1R_j^{(r-2)})\right)\\
&=&\sum_{i=1}^{b(B)}\sign(B)(-1)^{|B|+i-1}\cdot\Delta_B
\odot\left(\uuuu{p}_j+d^{-1}\paa C(\uuuu{c}_{k_1^B},...,
\widehat{\uuuu{c}_{k_i^B}},...,\uuuu{c}_{k_{b(B)}^B})\right).
\end{eqnarray*}
This is equal to 0, 
because for each pair $(i,l)\in\{1,...,b(B)\}^2$ with
$i<l$, the term
$\uuuu{p}_j+d^{-1} C(\uuuu{c}_{k_1^B},...,
\widehat{\uuuu{c}_{k_i^B}},...,\widehat{\uuuu{c}_{k_l^B}},
...,\uuuu{c}_{k_{b(B)}^B})$
and the term
$\uuuu{p}_{j+1}+d^{-1} C(\uuuu{c}_{k_1^B},...,
\widehat{\uuuu{c}_{k_i^B}},...,\widehat{\uuuu{c}_{k_l^B}},
...,\uuuu{c}_{k_{b(B)}^B})$
turn up twice and with different signs.
Therefore $\textup{Cont}(B,\paa_1R_j^{(r-2)})$ is a relative
cycle. In the proof of Theorem \ref{t3.3}, we found
\begin{eqnarray*}
&-&\textup{Cont}(B,\paa_1R_j^{(r-2)})
+\textup{Cont}(B,\paa_1R_{j+1}^{(r-2)}) \\
&=&\textup{Cont}(B,\paa_2R_j^{(r-1)})\\
&=&
\paa_2(\textup{Cont}(B,R_j^{(r-1)}))\\
&=&\paa (\textup{Cont}(B,R_j^{(r-1)}))
-\paa_1 (\textup{Cont}(B,R_j^{(r-1)})).
\end{eqnarray*}
In the last difference, the first term is a boundary,
and the second term is in $G_{n-r}$.
Therefore the class in $H_{n-2}(G_{n-r+1},G_{n-r})$ of
this difference is 0. Therefore the class in \eqref{4.51}
is independent of $j$. 

We see that all contributions to the class 
$[\paa_1 R_j^{(r-2)}]\in H_{n-2}(G_{n-r+1},G_{n-r})$ are
independent of $j$, so this class is independent of $j$. 
This implies the first statement in Claim 2. 

By \eqref{4.28} $0=\sum_{j=1}^d[X_j^{(r-2)}]\in E^1_{n,-1}$. Thus
\begin{eqnarray*}
0=\sum_{j=1}^d d^{r-1}_{n,-1}([X_j^{(r-2)}])=d\cdot 
[\paa_1 R_j^{(r-2)}]
\end{eqnarray*}
in the quotient in \eqref{4.49}. This is $\Leftarrow$ of
\eqref{4.50}. 

It rests to prove $\Rightarrow$ of \eqref{4.50}. 
Define $s:=n-r+2$. Then $r\in\{3,...,[\frac{n+3}{2}]\}$ implies
$s\in\{[\frac{n+2}{2}],...,n-1\}$. So part (a) and its proof 
apply. The map $\pr_\BB\circ d^1_{n-r+2,r-3}$ is injective,
and formula \eqref{4.43} gives its image in
$H_{n-2}(G_{n-r+1},G_{n-r})_\BB$. 
For $\Rightarrow$ of \eqref{4.50},
it is sufficient to show for $m\in\Z$
\begin{eqnarray}\label{4.52}
m\cdot \pr_\BB([\paa_1R_j^{(r-2)}]) \in 
\pr_\BB(\Imm(d^1_{n-r+2,r-3}))&\Rightarrow & d|m,
\end{eqnarray}
here $[\paa_1R_j^{(r-2)}]$ denotes the class in
$H_{n-2}(G_{n-r+1},G_{n-r})$ (not the class in the quotient
in \eqref{4.49}). 

For $\pr_\BB([\paa_1R_j^{(r-2)}])$, we need only the contributions 
of sets $B\in \BB(n-r+2)$.
By the proof of Theorem \ref{t3.3}, the contribution of such a 
set $B$ is 0 if $\beta_2(B)\neq n$, and it is given by \eqref{3.25}
if $\beta_2(B)=n$. We obtain
\begin{eqnarray*}
\pr_\BB([\paa_1R_j^{(r-2)}])
=\sum_{B\in\BB(n-r+2):\, \beta_2(B)=n}
(\pm 1)\uuuu{a}_1\uuuu{a}_{k_2^B}...\uuuu{a}_{k_{b(B)}^B}
\cdot \left[\Delta_B\odot T_B\right].
\end{eqnarray*}
Suppose that the left hand side of \eqref{4.52} holds
for some $m\in\Z$. 
Write $\lambda_B:=m\cdot (\pm 1)\uuuu{a}_1\uuuu{a}_{k_2^B}...
\uuuu{a}_{k_{b(B)}^B}$. Then we have a linear
combination $\sum_{A\in\AAA_1(n-r+2)}z_A\cdot[\Delta_A\odot T_A]$
with
\begin{eqnarray}\label{4.53}
m\cdot \pr_\BB([\paa_1R_j^{(r-2)}])
&=&\sum_{B\in \BB(n-r+2):\, \beta_2(B)=n} 
\lambda_B\cdot \left[\Delta_B\odot T_B\right]\\
&=&\pr_\BB\bigl(d^1_{n-r+2,r-3}\left(
\sum_{A\in\AAA_1(n-r+2)}z_A\cdot[\Delta_A\odot T_A]\right)\bigr).
\nonumber
\end{eqnarray}
Choose a set $B\in\BB(n-r+2)$ with $\beta_2(B)=n$. 
Claim 1 provides a sequence $(B_i)_{i\in N}$ of elements
$B_i\in\BB(n-r+2)$ with $B=B_n$, \eqref{4.44} and \eqref{4.45}.
Then \eqref{4.53} and \eqref{4.43} give the first line
of \eqref{4.46}, so 
\begin{eqnarray*}
z_{B^{(2)}}=\Bigl(\prod_{i=1}^{n-1}a_i\Bigr)\cdot z_{B^{(1)}}.
\end{eqnarray*}
By \eqref{4.43}, the coefficient $\lambda_B$ is 
\begin{eqnarray}
\lambda_B&=& \sign(B,B^{(1)})\cdot \Bigl(
z_{B^{(1)}}-(-1)^na_nz_{B^{(2)}}\Bigr)
\label{4.54}\\
&=& \sign(B,B^{(1)})\cdot z_{B^{(1)}}\cdot (1-\uuuu{a}_{n+1})
\nonumber\\
&=&\sign(B,B^{(1)})\cdot z_{B^{(1)}}\cdot (-1)^{n+1}d.
\nonumber
\end{eqnarray}
Therefore $\lambda_B\in d\Z$. Observe 
$\gcd(d, \uuuu{a}_1\uuuu{a}_{k_2^B}...\uuuu{a}_{k_{b(B)}^B})=1$.
Therefore, $d|m$. \eqref{4.52} is proved.
This finishes the proof of part (d).\hfill$\Box$

\begin{theorem}\label{t4.6}
The homology group $H_{n-1}(G,\Z)=H_{n-1}(G)$ is a $\Z$-lattice
of rank $\mu=(-1)^n\uuuu{a}_{n+1}=d+(-1)^n$. 
A $\Z$-basis of it is 
\begin{eqnarray}\label{4.55}
\begin{array}{ll}
[X_1^{(\frac{n}{2})}],...,[X_{d-1}^{(\frac{n}{2})}],
{}[\Delta_{A_{od}}\odot T_{A_{od}}], 
{}[\Delta_{A_{ev}}\odot T_{A_{ev}}],& 
\textup{if }n\textup{ is even,}\\
{}[X_1^{(\frac{n+1}{2})}],...,[X_{d-1}^{(\frac{n+1}{2})}],& 
\textup{if }n\textup{ is odd.}
\end{array}
\end{eqnarray}
Together with $[X_d^{([\frac{n+1}{2}])}]$, these elements satisfy
the relation \eqref{3.18} if $n$ is odd and the relation
\eqref{3.19} if $n$ is even. 
\end{theorem}

{\bf Proof:}
By Lemma \ref{t3.2}, in the case $n$ even, 
$[\Delta_{A_{od}}\odot T_{A_{od}}]$ and 
$[\Delta_{A_{ev}}\odot T_{A_{ev}}]$ are homology classes
in $H_{n-1}(G)$. By Theorem \ref{t3.3}, 
$[X_1^{([\frac{n+1}{2}])}],...,[X_d^{([\frac{n+1}{2}])}]$
are homology classes in $H_{n-1}(G)$, the relation
\eqref{3.18} holds if $n$ is odd, and the relation 
\eqref{3.19} holds if $n$ is even.
We have to show that the elements in \eqref{4.55} form a
$\Z$-basis of $H_{n-1}(G)$. 

First we consider even $n$. By Lemma \ref{t4.3} (b) and (c)
and Theorem \ref{t4.5} (a), 
$E^\infty_{s,n-1-s}=0$ for all $s$ except $s=n$
and $s=\frac{n}{2}$. 
Therefore $F^s_{n-1}=0$ for all $s<\frac{n}{2}$, 
and by Lemma \ref{t4.3} (c) 
\begin{eqnarray}
F^{\frac{n}{2}}_{n-1}&=&F^{\frac{n}{2}+1}_{n-1}=...=F^{n-1}_{n-1}
=\Imm(H_{n-1}(G_{\frac{n}{2}})\to H_{n-1}(G))\nonumber \\
&\cong & E^\infty_{\frac{n}{2},\frac{n}{2}-1}
=E^1_{\frac{n}{2},\frac{n}{2}-1}
=H_{n-1}(G_{\frac{n}{2}},G_{\frac{n}{2}-1})\cong\Z^2\label{4.56}
\end{eqnarray}
with generators the classes in $H_{n-1}(G)$ respectively in
$H_{n-1}(G_{\frac{n}{2}},G_{\frac{n}{2}-1})$ of the spheres
$\Delta_{A_{od}}\odot T_{A_{od}}\subset \C^{A_{od}}$ and 
$\Delta_{A_{ev}}\odot T_{A_{ev}}\subset \C^{A_{ev}}$. 
By Theorem \ref{t4.5} (d) and Lemma \ref{t4.4} 
\begin{eqnarray}\label{4.57}
E^\infty_{n,-1} =\www{E}^{\frac{n+2}{2}}_{n,-1}=
\bigoplus_{j=1}^{d-1}\Z\cdot [X_j^{(\frac{n}{2})}]
\subset E^1_{n,-1}=H_{n-1}(G,G_{n-1}),
\end{eqnarray}
$E^\infty_{n,-1}$ is a $\Z$-lattice of rank $d-1$, and
\begin{eqnarray}\label{4.58} 
E^\infty_{n,-1}\cong \frac{F^n_{n-1}}{F^{\frac{n}{2}}_{n-1}}
=\frac{H_{n-1}(G)}{\Imm(H_{n-1}(G_{\frac{n}{2}})\to H_{n-1}(G)))}.
\end{eqnarray}
Therefore $H_{n-1}(G)$ is a $\Z$-lattice of rank $d+1$
with $\Z$-basis as in \eqref{4.55}, and 
$F^{\frac{n}{2}}_{n-1}$ is a primitive sublattice
of rank 2 with generators the classes of the spheres 
$\Delta_{A_{od}}\odot T_{A_{od}}$ and 
$\Delta_{A_{ev}}\odot T_{A_{ev}}$. 

Now we consider odd $n$. By Lemma \ref{t4.3} (b) and Theorem
\ref{t4.5} (a), $E^\infty_{s,n-1-s}=0$ for all $s\neq n$. 
Therefore $F^s_{n-1}=0$ for all $s\neq n$ and
\begin{eqnarray}\label{4.59}
H_{n-1}(G)=F^n_{n-1}\cong E^\infty_{n,-1}.
\end{eqnarray}
By Theorem \ref{t4.5} (d) and Lemma \ref{t4.4} 
\begin{eqnarray}\label{4.60}
E^\infty_{n,-1} =\www{E}^{\frac{n+3}{2}}_{n,-1}=
\bigoplus_{j=1}^{d-1}\Z\cdot [X_j^{(\frac{n+1}{2})}]
\subset E^1_{n,-1}=H_{n-1}(G,G_{n-1}),
\end{eqnarray}
and $E^\infty_{n,-1}$ is a $\Z$-lattice of rank $d-1$.
Therefore $H_{n-1}(G)$ is a $\Z$-lattice of rank $d-1$
with $\Z$-basis as in \eqref{4.55}.
\hfill$\Box$

\begin{remarks}\label{t4.7}
In the proof of part (d) of Theorem \ref{t4.5},
we considered the contribution 
$\textup{Cont}(B,\paa_1R_j^{(r-2)})$
of a set $B\in \AAA_1(n-r+1)$ ($B$ is thick with 
$|B|=n-r+1$) to the chain $\paa_1R_j^{(r-2)}$. We showed
that it is a relative cycle in $H_{n-2}(G_{n-r+1},G_{n-r})$,
and we saw that it is independent of $j$. 
We did not make precise which cycle it is, as we did not need 
this. Now we will say which cycle it is, but leave the proof
to the reader:
\begin{eqnarray}\label{4.61}
&&[\textup{Cont}(B,\paa_1R_j^{(r-2)})]\\
&=&\varepsilon_1\cdot [\Delta_B\odot \uuuu{e}\bigl(
d^{-1}C(\uuuu{c}_{k_2^B}-\uuuu{c}_{k_1^B},...,
\uuuu{c}_{k_{b(B)}^B}-\uuuu{c}_{k_{b(B)-1}^B})\bigr)]\nonumber\\
&=&\varepsilon_2\cdot [\Delta_B\odot\uuuu{e}\bigl( C(
\uuuu{b}_{k_2^B}-\uuuu{b}_{k_1^B},...,
\uuuu{b}_{k_{b(B)}^B}-\uuuu{b}_{k_{b(B)-1}^B})\bigr)].\nonumber\\
&=&\varepsilon_2\cdot 
\uuuu{a}_{k_1^B}\uuuu{a}_{k_2^B}...\uuuu{a}_{k_{b(B)-1}^B}\cdot 
[\Delta_B\odot\uuuu{e}\bigl( C(
\uuuu{d}_1^B,..., \uuuu{d}_{b(B)-1}^B)\bigr)].\nonumber\\
&&\textup{with }\varepsilon_1=\sign(B)(-1)^{|B|},\ 
\varepsilon_2=\varepsilon_1 \cdot(-1)^{n(b(B)-1)}.\nonumber
\end{eqnarray}
\end{remarks}

\section{Proof of the main result}\label{c5}
\setcounter{equation}{0}

\noindent
Throughout this section we fix 
a cycle type singularity $f(x_1,...,x_n)$ as in
\eqref{1.2} with $n\geq 2$ and $a_1,...,a_n\in\N$ with
\eqref{1.3}, as in the sections \ref{c3} and \ref{c4}.
We want to prove the main result of this paper, 
Theorem \ref{t1.3}. 

The deformation lemma in \cite[3.]{Co82}
implies that the inclusion map $i_g:G\hookrightarrow \oooo{F_0}$
induces an epimorphism 
$(i_g)_*:H_{n-1}(G,\Z)\to H_{n-1}(\oooo{F_0},\Z)$.
Both groups are $\Z$-lattices of rank $\mu$. This holds
by Theorem \ref{t4.6} for $H_{n-1}(G,\Z)$ and by
\cite{Mi68} for $H_{n-1}(\oooo{F_0},\Z)$.
Therefore the epimorphism is an
\begin{eqnarray*}
\textup{isomorphism }\quad
(i_g)_*:H_{n-1}(G,\Z)\to H_{n-1}(\oooo{F_0},\Z).
\end{eqnarray*}

It rests to determine the monodromy on $H_{n-1}(G,\Z)$.
We call this monodromy $h_{mon}$. 
First we give the weights
$(w_1,...,w_n)$ of the quasihomogeneous polynomial $f$
(see e.g. Lemma 4.1 in \cite{HZ19}):
\begin{eqnarray}
w_j&=&\frac{v_j}{d}\in\Q\cap (0,1)
\quad\textup{for }j\in N\quad \textup{ with}\nonumber\\
d&=&\prod_{j=1}^na_j-(-1)^n=\mu-(-1)^n,\nonumber\\
v_j&=&\sum_{l=1}^n(-1)^{l-1}
\prod_{k=j+l}^{j+n-1}a_{(k)_{\mmod n}},
\textup{ with}\label{5.1}\\
&&a_jv_j+v_{(j+1)_{\mmod n}}=d.\nonumber
\end{eqnarray}
One sees easily $\gcd(v_1,d)=...=\gcd(v_n,d)$. 
We defined already $b:=d/\gcd(v_1,d)$. 
In the case of a quasihomogeneous singularity, one can
give explicitly a diffeomorphism $\Phi_{mon}:\C^n\to\C^n$ 
of the Milnor fibration which induces the monodromy 
on the homology of the Milnor fiber and on $H_{n-1}(G,\Z)$. 
It looks as follows.
\begin{eqnarray*}
\Phi_{mon}:\C^n\to\C^n,\quad (z_1,..,z_n)\mapsto
(e^{2\pi i w_1}z_1,...,e^{2\pi i w_n}z_n).
\end{eqnarray*}
It maps the point $\uuuu{e}(\uuuu{p}_j)$ to the point
$\uuuu{e}(\uuuu{p}_{(j+v_1)_{\mmod d}})$, because of
$a_jv_j+v_{(j+1)_{\mmod n}}=d$. 
Therefore it maps the chain $R_j^{(k)}$ to the chain 
$R_{(j+v_1)_{\mmod d}}^{(k)}$
and the chain $X_j^{(k)}$ to the chain 
$X_{(j+v_1)_{\mmod d}}^{(k)}$.
In the case of even $n$, it maps the spheres
$\Delta_{A_{od}}\odot T_{A_{od}}\subset \C^{A_{od}}$ and 
$\Delta_{A_{ev}}\odot T_{A_{ev}}\subset \C^{A_{ev}}$ 
to themselves. 

Now we can apply the results of section \ref{c2}.
In order to make the relation to these results transparent, 
we introduce the following notations,
\begin{eqnarray*}
\delta_j&:=& [X_j^{([\frac{n+1}{2}])}]\in H_{n-1}(G,\Z)
\quad\textup{for }j\in\{1,...,d\},
\end{eqnarray*}
for odd $n$ or even $n$, and
\begin{eqnarray*}
\gamma&:=& [\Delta_{A_{od}}\odot T_{A_{od}}] +
(-1)^{\frac{n}{2}}a_1a_3...a_{n-1}\cdot 
[\Delta_{A_{ev}}\odot T_{A_{ev}}]\in H_{n-1}(G,\Z),\\
\beta&:=& [\Delta_{A_{ev}}\odot T_{A_{ev}}]\in H_{n-1}(G,\Z),\\
c&:=& (-1)^{(\frac{n}{2}+2)(\frac{n}{2}+1)\frac{1}{2}}\cdot 
\uuuu{a}_1\uuuu{a}_3...\uuuu{a}_{n-1}\in\Z,
\end{eqnarray*}
for even $n$. Then
\begin{eqnarray*}
H_{n-1}(G,\Z)&=&\bigoplus_{j=1}^{d-1}\Z\cdot\delta_j
\oplus\left\{\begin{array}{ll}
0&\textup{if }n\textup{ is odd},\\
\Z\cdot\gamma\oplus\Z\cdot\beta& 
\textup{if }n\textup{ is even},\end{array}\right.\\
h_{mon}(\delta_j)&=&\delta_{(j+v_1)_{\mmod d}},\\
h_{mon}(\gamma)&=&\gamma,\quad h_{mon}(\beta)=\beta,
\quad\textup{if }n\textup{ is even,}\\
\sum_{j=1}^d\delta_j&=& \left\{
\begin{array}{ll}
0&\textup{if }n\textup{ is odd (because of \eqref{3.18}),}\\
c\cdot \gamma&\textup{if }n\textup{ is even
(because of \eqref{3.19}).}
\end{array}\right. 
\end{eqnarray*}
The proof of Lemma \ref{t2.3} can be adapted to show
for odd $n$
\begin{eqnarray*}
(H_{n-1}(G,\Z)),h_{mon})\cong 
(\gcd(v_1,d)-1)\Or(t^b-1)\oplus\Or(\frac{t^b-1}{t-1}).
\end{eqnarray*}
For even $n$, we see with Definition \ref{t2.1},
Lemma \ref{t2.3} and Lemma \ref{t2.4}
\begin{eqnarray*}
&&(H_{n-1}(G,\Z),h_{mon})\\
&\stackrel{\textup{Definition \ref{t2.1}}}{\cong}&
(H^{(d,c)},(h^{(d,c)})^{v_1}) \oplus\Or(t-1)\\
&\stackrel{\textup{Lemma \ref{t2.3}}}{\cong}&
(\gcd(v_1,d)-1)\Or(t^b-1)\oplus \textup{Lo}^{(b,c)}
\oplus \Or(t-1)\\
&\stackrel{\textup{Lemma \ref{t2.4}}}{\cong}&
\gcd(v_1,d)\Or(t^b-1)\oplus\Or(t-1),
\end{eqnarray*}
the last isomorphism uses $\gcd(b,c)=1$ and \eqref{2.3}.
This finishes the proof of Theorem \ref{t1.3}.

\begin{remarks}\label{t5.1}
Cooper's paper \cite{Co82} studied the integral monodromy 
of the cycle type singularities. It is split into 22
sections. The sections 1 und 2 are an introduction
and a discussion of a degenerate case.
The sections 3 to 14 are devoted to the deformation lemma
in section 3, which yields that $(i_g)_*$ is surjective.

Section 15 introduces the deformed simplices $\Delta_A$,
the tori $T_A$, the chains $\Delta_A\odot T_A$,
the blocks of sets $A\subset N$, and it gives most of
the statements of Lemma \ref{t3.2}.

Section 16 states formulas which are close to our
formula \eqref{4.40}. 

At the end of section 17,
a formula which is close to our formula \eqref{4.43}
is derived. Section 17 also defines the subsets $G_s$
and thick sets $A$ and states a part of the formulas
\eqref{4.7}--\eqref{4.13} for the relative homology groups.

Section 18 first introduces the spectral sequence 
$(E^r_{s,t},d^r_{s,t})$ of the filtration 
$G=G_n\supset G_{n-1}\supset ...\supset G_1\supset G_0=\emptyset$.
But it gives neither a reference, nor the properties
in Theorem \ref{t4.2}.
Then it states and proves the injectivity of the
maps $d^1_{s,n-1-s}$ in our Theorem \ref{t4.5} (a). 
Though it does not use our Claim 1, but uses
a specific sequence of $n\cdot (n-s)$ elements 
$B_i\in\BB(s)$.
Claim 1 is more efficient and can be used also for the
proof of Theorem \ref{t4.5} (d).

The lemma in Section 19 is crucial, but not correct.
It claims correctly $[X_j^{(1)}]\in E^2_{n,-1}$,
but wrongly $[X_j^{(2)}]\in E^\infty_{n,-1}$.
By Theorem \ref{t4.5} (d), this is true only if 
$n\in\{3,4\}$. 
The wrong argument is in the proof of the lemma at the
end of section 19: A part of the boundary is missed.
This makes Cooper believe that $X_j^{(2)}$ is always a cycle.
He constructs and sees only the beginning for
$k\in\{0,1,2\}$ of the sequence 
$(X_j^{(k)})_{k\in\{1,...,[\frac{n+1}{2}]\}}$
of chains. 

The lemma in section 20 has some similarity with Claim 2
(in the proof of Theorem \ref{t4.5} (d)). But it is too
special. 

The lemma in section 21 concludes wrongly that
the classes $[X_j^{(2)}]_{j=1,...,d-1}$ form a
$\Z$-basis of $E^\infty_{n,-1}$. 
This is true only for $n\in\{3,4\}$. 
It builds on the lemma
in section 19. The wrong statements on the elements
and generators of $E^\infty_{n,-1}$ in the sections 19
and 21 form the first serious mistake.

The final section 22 makes the second serious mistake,
in the case of even $n$. It states correctly that
for even $n$ $E^\infty_{s,n-1-s}\neq 0$ only if 
$s\in\{\frac{n}{2},n\}$. And it derives correctly
$(E^\infty_{\frac{n}{2},\frac{n}{2}-1},h_{mon})
\cong 2\Or(t-1)$ and 
\begin{eqnarray}\label{5.1}
(E^\infty_{n,-1},h_{mon})\cong (\gcd(v_1,d)-1)\Or(t^b-1)
\oplus \Or(\frac{t^b-1}{t-1}).
\end{eqnarray}
But then it supposes wrongly that the quotients 
$E^\infty_{\frac{n}{2},\frac{n}{2}-1}$ and 
$E^\infty_{n,-1}$ of the filtration $F^\bullet_{n-1}$ 
of $H_{n-1}(G)$ lift to a splitting of $H_{n-1}(G)$ 
which is invariant by the monodromy.
This would give a non-standard decomposition of 
$(H_{n-1}(G),h_{mon})$ into Orlik blocks 
(and this is what Cooper claims to have for even $n$). 

But Cooper's paper contains most of the ingredients which
we need for the proof of Orlik's conjecture. 
We corrected the two mistakes, we were more careful with
the signs and proofs, we introduced more notations, we had
the full sequence of chains 
$(X_j^{(k)})_{k\in\{1,...,[\frac{n+1}{2}]\}}$,
and we had the algebraic statements in our section \ref{c2}.
But Cooper's paper \cite{Co82} was an indispensable basis on which
we could build.
\end{remarks}

\begin{remarks}\label{t5.2}
(i) Orlik and Randell study in \cite{OR77} mainly the 
integral monodromy of the chain type singularities.
But section 3 says also something about the cycle type
singularities. In (3.4.1) they make a conjecture on the
integral monodromy of the cycle type singularities which
would imply Orlik's conjecture. 
Their ansatz for a proof is very different from \cite{Co82}.
They point themselves to a gap in this ansatz.

(ii) \cite{HM20} combines the result in \cite{OR77}
for the chain type singularities with an algebraic result.
Together they imply that Orlik's conjecture holds for the
chain type singularities. And \cite{HM20} gives a relative
result, which says that the Thom-Sebastiani sum of any
two quasihomogeneous singularities satisfies Orlik's 
conjecture if the two singularities satisfy Orlik's conjecture.
As Orlik's conjecture is valid for the chain type singularities
by \cite{HM20} and for the cycle types singularities by this
paper, it is valid also for all their (iterated) 
Thom-Sebastiani sums.
\end{remarks}

\end{document}